\documentclass{article}
\usepackage{amssymb}
\usepackage{amsfonts}
\usepackage{amsmath}

\newtheorem{theorem}{Theorem}

\newtheorem{corollary}[theorem]{Corollary}

\newtheorem{lemma}[theorem]{Lemma}

\begin{document}

\title{Macroscopic and microscopic structures of the family tree for
decomposable critical branching processes}
\author{Vladimir Vatutin\thanks{
Steklov Mathematical institute RAS, Gubkin str. 8, Moscow, 119991,
Russia; e-mail: vatutin@mi.ras.ru}}
\date{}
\maketitle

\begin{abstract}
A decomposable strongly critical Galton-Watson branching process with $N$
types of particles labelled $1,2,...,N$ is considered in which a type~$i$
parent may produce individuals of types $j\geq i$ only. This model may be
viewed as a stochastic model for the sizes of a geographically structured
population occupying $N$ islands, the location of a particle being
considered as its type. The newborn particles of island $i\leq N-1$ either
stay at the same island or migrate, just after their birth to the islands $%
i+1,i+2,...,N$. Particles of island $N$ do not migrate. We investigate the
structure of the family tree for this process, the distributions of the
birth moment and the type of the most recent common ancestor of the
individuals existing in the population at a distant moment $n.$
\end{abstract}

\section{Introduction and main results}

We consider a Galton-Watson branching process with $N$ types of particles
labelled $1,2,...,N$ \ and denote by
\begin{equation*}
\mathbf{Z}(n)=(Z_{1}(n),...,Z_{N}(n)),\quad \mathbf{Z}(0)=(1,0,...,0)
\end{equation*}%
the population vector at time $n\in \mathbb{Z}_{+}=\left\{ 0,1,...\right\} $%
. Along with $\mathbf{Z}(n)$ we deal with the process \
\begin{equation*}
\mathbf{Z}(m,n)=(Z_{1}(m,n),...,Z_{N}(m,n)),
\end{equation*}%
where $Z_{i}(m,n)$ is the number of type $i$ particles existing in $\mathbf{Z%
}(\cdot )$ at moment $m<n$ and having nonempty number of descendants at
moment $n$. We agree to write $Z_{i}(n,n)=Z_{i}(n)$.

The process $\mathbf{Z}( \cdot ,n) $ is called a reduced branching process
and can be thought of as the family tree relating the individuals alive at
time $n$. An important characteristic of the reduced process is the birth
moment $\beta _{n}$ of the most recent common ancestor (MRCA) of all
individuals existing in the population at moment $n$ defined as \
\begin{equation*}
\beta _{n}=\max \left\{ m\leq n-1:Z_{1}( m,n) +Z_{2}( m,n)
+...+Z_{N}(m,n)=1\right\} .
\end{equation*}

The structure of the family tree and the asymptotic distribution of the
birth moment of the MRCA for single-type Galton-Watson branching processes
have been studied in \cite{FP},\cite{FZ},\cite{LS} and \cite{Zub}. The case
of multitype indecomposable critical Markov branching processes was
considered in \cite{Yak83}. Family trees for more general models of
branching processes were investigated in \cite{BorVat}, \cite{FVat99},\cite%
{Sag85},\cite{Sag87},\cite{Sag95},\cite{Vat2003},\cite{VD06},\cite{VV79}.
However, the reduced processes for decomposable branching processes have not
been analyzed yet. We fill this gap in the present paper and study various
properties of the family tree for a particular case of the decomposable
Galton-Watson branching processes. Namely, we consider the Galton-Watson
branching process with $N$ types of particles labelled $1,2,...,N$ in which
a type $i$ parent may produce individuals of types $j\geq i$ only. This
model may be viewed as a stochastic model for the sizes of a geographically
structured population occupying $N$ islands, the location of a particle
being considered as its type. The reproduction laws of particles depend on
the island on which the particles are located. The newborn particles of
island $i\leq N-1$ either stay at the same island or migrate, just after
their birth to the islands $i+1,i+2,...,N$. Particles of island $N$ do not
migrate.

We investigate the structure of the family tree of this process, the
distributions of the birth moment $\beta _{n}$ and the type $\zeta _{n}$ of
the MRCA. It is shown, in particular, that, as $n\rightarrow \infty $ the
conditional reduced process
\begin{equation*}
\left\{ \mathbf{Z}(n^{t}\log n,n),0\leq t<1|\mathbf{Z}(n)\neq \mathbf{0}%
\right\}
\end{equation*}%
converges in a certain sense to an $N-$dimensional inhomogeneous branching
process $\left\{ \mathbf{R}(t),0\leq t<1\right\} $ which, for $t\in \lbrack
0,2^{-(N-1)})$ consists of a single particle of type $1$ only and for $t\in
\lbrack 2^{-(N-i+1)},2^{-(N-i)}),i=2,...,N$ consists of type $i$ particles
only. These particles are born at moment $t=2^{-(N-i+1)}$ and die at moment $%
t=2^{-(N-i)}$ producing at this moment a random number of descendants having
type $\min (i+1,N)$. This gives a macroscopic view on the structure of the
family tree of the process.

On the other hand, for each $i=1,2,...,N-1$ the conditional process
\begin{equation*}
\left\{ \mathbf{Z}((y+(\log n)^{-1})n^{2^{-(N-i)}},n),0<y<\infty \big|\,%
\mathbf{Z}(n)\neq \mathbf{0}\,\right\}
\end{equation*}%
converges in a certain sense, as $n\rightarrow \infty $ to a continuous-time
homogeneous Markov branching process $\left\{ \mathbf{U}_{i}(y),0\leq
y<\infty \right\} $ which is initiated at time $y=0$ by a random number of
type $i$ particles. These type $i$ particles have an exponential life-length
distribution. Dying each of them produces either two particles of type $i$
or one particle of type $i+1$ (both options with probability 1/2). Particles
of type $i+1$ in this process are immortal and produce no offspring. This
provides a microscopic view on the structure of the family tree.

To present our results in a more formal way we need some notation. Let $%
\mathbf{e}_{i}$ be a vector whose $i$-th component is equal to one while the
remaining are zeros. The first moments of the components of $\mathbf{Z}(n)$
will be denoted as%
\begin{equation*}
m_{ij}(n)=\mathbf{E}\left[ Z_{j}(n)|\mathbf{Z}\left( 0\right) =\mathbf{e}_{i}%
\right]
\end{equation*}%
with $m_{ij}=m_{ij}(1)$ \thinspace being the average number of children of
type $j$ produced by a particle of type $i$.

Since $m_{ij}=0$ if $i>j$, the mean matrix $\mathbf{M}$ of the decomposable
Galton-Watson branching process has the form
\begin{equation}
\mathbf{M=}(m_{ij})_{i,j=1}^{N}=\left(
\begin{array}{ccccc}
m_{11} & m_{12} & ... & ... & m_{1N} \\
0 & m_{22} & ... & ... & m_{2N} \\
0 & 0 & m_{33} & ... & ... \\
... & ... & ... & ... & ... \\
... & ... & ... & ... & ... \\
0 & 0 & ... & 0 & m_{NN}%
\end{array}%
\right) .  \label{matrix1}
\end{equation}

To go further it is convenient to deal with the probability generating
functions for the reproduction laws of particles

\begin{equation}
h_{i}(s_{1},...,s_{N})=\mathbf{E}\left[ s_{i}^{\eta _{ii}}...\,s_{N}^{\eta
_{iN}}\right] ,\ i=1,2,...,N,  \label{DefNONimmigr}
\end{equation}%
where $\eta _{ij}$ represent the numbers of daughters of type $j$ a mother
of type $i$.

We say that \textbf{Hypothesis A} is valid if the $N-$type decomposable
process is strongly critical, i.e. (see \cite{FN2}),
\begin{equation}
m_{ii}=\mathbf{E}\left[ \eta _{ii}\right] =1,\ i=1,2,...,N,  \label{Matpos}
\end{equation}%
and, in addition,%
\begin{equation}
m_{i,i+1}=\mathbf{E}\left[ \eta _{i,i+1}\right] \in ( 0,\infty ) ,\
i=1,2,...,N-1,  \label{Maseq}
\end{equation}%
and%
\begin{equation}
\mathbf{E}\left[ \eta _{ij}\eta _{ik}\right] <\infty ,\,i=1,...,N;\
k,j=i,i+1,...,N  \label{FinCovar}
\end{equation}%
with%
\begin{equation}
b_{i}=\frac{1}{2}Var\left[ \eta _{ii}\right] \in ( 0,\infty ) ,\ i=1,2,...,N.
\label{FinVar}
\end{equation}%
Thus, a particle of the process under consideration is able to produce the
direct descendants of its own type, of the next in the order type, and (not
necessarily, as direct descendants) of all the remaining in the order types,
but not any preceding ones.

To simplify the presentation we fix, from now on $N\geq 2$ and use, when it
is convenient the notation%
\begin{equation*}
\gamma _{0}=0,\ \gamma _{i}=\gamma _{i}(N)=2^{-(N-i)},\ i=1,2,...,N.
\end{equation*}
We also suppose (if otherwise is not stated) that $\mathbf{Z}(0)=\mathbf{e}%
_{1}$, i.e., assume that the branching process under consideration is
initiated at time zero by a single particle of type $1$.

Let $\xi ^{(i)}(j),i=1,2,...,N;j=1,2,...$ be a tuple of independent
identically distributed random variables with probability generating function%
\begin{equation*}
f(s)=\mathbf{E}\left[ s^{\xi ^{(i)}(j)}\right] =1-\sqrt{1-s}.
\end{equation*}%
By means of the tuple we give a detailed construction of an $N-$type
decomposable branching process $\mathbf{R}(t)=(R_{1}(t),...,R_{N}(t)),0\leq
t<1,$ where $R_{i}(t)$ is the number of type $i$ individuals in the
population at moment~$t$. It is this process describes the macroscopic
structure of the family tree $\left\{ \mathbf{Z}(m,n),0\leq m\leq n\right\}$
as $n\rightarrow \infty $.

Let $\mathbf{R}(t)=\mathbf{e}_{1}$ for $\gamma _{0}\leq t<\gamma _{1}$
meaning that the branching process $\mathbf{R}(t)$ starts at $t=0$ by a
single individual of type $1$ which survives up to (but not at) moment $%
\gamma _{1}$ without reproduction. If $\gamma _{i}\leq t<\gamma
_{i+1},i=1,2,...,N-1$ then%
\begin{equation*}
R_{k}(t)=\left\{
\begin{array}{ccc}
\sum_{j=1}^{R_{i}(\gamma _{i}-0)}\xi ^{(i)}(j) & \text{if} & k=i+1 \\
&  &  \\
0 & \text{if} & k\neq i+1%
\end{array}%
\right. .
\end{equation*}%
Thus, within the interval $\gamma _{i}\leq t<\gamma _{i+1}$ the population
consists of type $i+1$ particles only. These particles were born at moment $%
\gamma _{i}-0$ by particles of type $i$ evolving without reproduction within
the interval $\gamma _{i-1}\leq t<\gamma _{i}$. More precisely, the $j-$th
particle of type $i$ produces at its death moment $\gamma _{i}-0$ a random
number $\xi ^{(i)}(j)$ children of type $i+1$ and no particles of other
types.

In what follows we use the symbol $\Longrightarrow $ to denote convergence
in the space $D_{[a,b)}(\mathbb{Z}_{+}^{N})$ of cadlag functions $\mathbf{x}%
(t),a\leq t<b$ with values in \ $\mathbb{Z}_{+}^{N}$ endowed with the metric
of Skorokhod topology. Besides, we agree to consider $\mathbf{Z}(x,n)$ as $%
\mathbf{Z}([x],n),$ where $[x]$ is the integer part of $x$.

For $0\leq t\leq 1$ put%
\begin{equation*}
g_{n}(t)=1_{\{0\leq t<\gamma _{1}\}}+g_{n}1_{\{\gamma _{1}\leq t\leq 1\}}
\end{equation*}%
where $g_{n}$ is a positive monotone increasing sequence such that
\begin{equation*}
\lim_{n\rightarrow \infty }g_{n}=\infty \text{ and }\lim_{n\rightarrow
\infty }n^{-\varepsilon }g_{n}=0\text{ for any }\varepsilon >0.
\end{equation*}

\begin{theorem}
\label{T_SkorohConst}Let Hypothesis A be valid. Then, as $n\rightarrow
\infty $

1) the finite-dimensional distributions of the process
\begin{equation*}
\left\{ (\mathbf{Z}(n^{t}g_{n}(t),n),0\leq t<1)|\mathbf{Z}(n)\neq \mathbf{0}%
\right\}
\end{equation*}%
converge to the finite-dimensional dis\-tri\-bu\-tions of $\left\{ \mathbf{R}%
(t),0\leq t<1\right\};$

2) for any $i=0,1,2,...,N-1$
\begin{equation*}
\mathcal{L}\left\{ (\mathbf{Z}(n^{t}g_{n}(t),n),\gamma _{i}\leq t<\gamma
_{i+1})\,|\,\mathbf{Z}(n)\neq \mathbf{0}\right\} \Longrightarrow \mathcal{L}%
\left\{ \mathbf{R}(t),\gamma _{i}\leq t<\gamma _{i+1}\right\} .
\end{equation*}
\end{theorem}

\textbf{Remark 1.} Theorem \ref{T_SkorohConst} shows that the passage to
limit under the macroscopic time-scaling $n^{t}g_{n}(t)$ transforms the
reduced process into an inhomogeneous branching process which consists at
any given moment of particles of a single type only. In particular, the
phase transition from type $i$ to type $i+1$ in the prelimiting process
happens, roughly speaking, at moment $n^{\gamma _{i}}$. This gives a
macroscopic view on the family tree of the reduced process. The microscopic
structure of the family tree described by Theorem \ref{T_Skhod1} below
clarifies the nature of the revealed phase transition.

Let $c_{ji},1\leq j\leq i\leq N$ be a tuple of positive numbers in which $%
c_{ii}=b_{i}^{-1}$ for $i=1,2,...,N$ and%
\begin{equation}
c_{ji}=\sqrt{b_{j}^{-1}m_{j,j+1}c_{j+1,i}}\text{ for }\ j\leq i-1,\quad
C_{i}=c_{1i}.  \label{Const1}
\end{equation}

It is not difficult to check that%
\begin{equation}
c_{iN}=\left( \frac{1}{b_{N}}\right) ^{1/2^{N-i}}\prod_{j=i}^{N-1}\left(
\frac{m_{j,j+1}}{b_{j}}\right) ^{1/2^{j-i+1}}.  \label{Const2}
\end{equation}

We now define a tuple of continuous time Markov processes%
\begin{eqnarray*}
\mathbf{U}_{i}(y) &=&(U_{i1}(y),...,U_{iN}(y)),\ 0\leq y<\infty ,\
i=1,2,...,N-1, \\
&& \\
\mathbf{U}_{N}(x) &=&(U_{N1}(x),...,U_{NN}(x)),\ 0\leq x<1.
\end{eqnarray*}

First we describe the structure of the processes $\mathbf{U}_{i}(y),1\leq
i\leq N-1$. In this case $U_{ij}(y)\equiv 0,~0\leq y<\infty ,~j\neq i,i+1,$
while the pair
\begin{equation*}
(U_{ii}(y),U_{i,i+1}(y)),0\leq y<\infty ,
\end{equation*}%
constitutes a two-type continuous-time homogeneous Markov branching process
with particles of types $i$ and $i+1$. This two-type process is initiated at
time $y=0$ by a random number $R_{i}$ of type $i$ particles whose
distribution is specified by the probability generating function%
\begin{equation}
\mathbf{E}\left[ s_{i}^{R_{i}}\right] =\mathbf{E}\left[ s_{i}^{U_{ii}(0)}%
\right] =1-(1-s_{i})^{1/2^{i-1}}  \label{Dist_ro}
\end{equation}%
(in particular, $U_{11}(0)=1$ with probability 1). The life-length
distribution of type $i$ particles is exponential with parameter $%
2b_{i}c_{iN}$. Dying each particle of type $i$ produces either two particles
of its own type or one particle of type $i+1$ (each option with probability
1/2). Particles of type $i+1$ of $\mathbf{U}_{i}(\cdot )$ are immortal and
produce no children.

The structure of the $N-$ dimensional process $\mathbf{U}_{N}(x),0\leq x<1$
is different. If $j<N$ then $U_{Nj}(x)\equiv 0,~0\leq x<1,$ while the
component $U_{NN}(\cdot )$ is a single-type inhomogeneous Markov branching
process initiated at time $x=0$ by a random number $R_{N}$ of type $N$
individuals distributed in accordance with probability generating function

\begin{equation}
\mathbf{E}\left[ s_{N}^{R_{N}}\right] =\mathbf{E}\left[ s_{N}^{U_{NN}(0)}%
\right] =1-(1-s_{N})^{1/2^{N-1}}.  \label{DefroN}
\end{equation}%
The life-length of each of $R_{N}$ type $N$ initial particles is uniformly
distributed on the interval $\left[ 0,1\right] $. Dying such a particle
produces exactly two children of type~$N$ and nothing else. If the death
moment of the parent particle is $x$ then the life length of each of its
offspring has the uniform distribution on the interval $[x,1]$
(independently of the behavior of other particles and the prehistory of the
process). Dying each particle of the process produces exactly two
individuals of type $N$ and so on...\thinspace .

We are now ready to formulate one more important result of the paper,
describing the microscopic structure of the family tree.

Let $l_{n}$ be a monotone decreasing sequence such that
\begin{equation*}
\lim_{n\rightarrow \infty }l_{n}=0\text{ and }\lim_{n\rightarrow \infty
}n^{\varepsilon }l_{n}=\infty \text{ for any }\varepsilon >0.
\end{equation*}

\begin{theorem}
\label{T_Skhod1}Let Hypothesis A be valid. Then, as $n\rightarrow \infty $

1) for each $i=1,2,...,N-1$
\begin{equation*}
\mathcal{L}\left\{ (\mathbf{Z}\left( (y+l_{n})n^{\gamma _{i}},n\right)
,0\leq y<\infty )\big|\,\mathbf{Z}(n)\neq \mathbf{0}\right\} \Longrightarrow
\mathcal{L}_{R_{i}}\left\{ \mathbf{U}_{i}(y),0\leq y<\infty \,\right\} ,
\end{equation*}%
where $\mathcal{L}_{R_{i}}$ means that $\mathbf{U}_{i}(\cdot )$ is initiated
at time $y=0$ by a random number $R_{i}$ particles of type $i$ (with $%
R_{1}\equiv 1)$;
\begin{equation*}
\mathit{2)}\,\mathcal{L}\left\{(\mathbf{Z}((x+l_{n})n,n),0\leq x<1)\,|\,%
\mathbf{Z}(n)\neq \mathbf{0}\right\} \Longrightarrow \mathcal{L}%
_{R_{N}}\left\{ \mathbf{U}_{N}(x),0\leq x<1\,\right\},
\end{equation*}%
where $\mathcal{L}_{R_{N}}$ means that $\mathbf{U}_{N}(\cdot)$ is initiated
at time $x=0$ by a random number $R_{N}$ particles of type $N.$
\end{theorem}

\textbf{Remark 2.} Theorems \ref{T_SkorohConst} and \ref{T_Skhod1} reveal an
interesting phenomenon in the development of the critical decomposable
branching processes which may be expressed in terms of the "island"
interpretation of the processes as follows: If the population survives up to
a distant moment $n$, then all surviving individuals are located at this
moment on island $N$ and, moreover, at each moment in the past their
ancestors were (asymptotically) located not more than on two specific
islands.

Basing on the conclusions of Theorems \ref{T_SkorohConst} and \ref{T_Skhod1}
we give in the next theorem an answer to the following important question:
what is the asymptotic distribution of the birth moment of the MRCA for the
population survived up to a distant moment $n$?

\begin{theorem}
\label{T_mrcaMany}Let Hypothesis A be valid. Then

1)
\begin{equation*}
\lim_{n\rightarrow \infty }\mathbf{P}\left( \beta _{n}\ll n^{\gamma _{1}}%
\big|\,\mathbf{Z}(n)\neq \mathbf{0}\right) =0;
\end{equation*}%
2) if $y\in (0,\infty )$ then for $i=1,2,...,N-1$
\begin{equation*}
\lim_{n\rightarrow \infty }\mathbf{P}\left( \beta _{n}\leq yn^{\gamma _{i}}%
\big|\,\mathbf{Z}(n)\neq \mathbf{0}\right) =1-\frac{1}{2^{i}}-\frac{1}{2^{i}}%
e^{-2b_{i}c_{iN}y};
\end{equation*}%
3) for $i=1,2,...,N-1$
\begin{equation}
\lim_{n\rightarrow \infty }\mathbf{P}\left( \beta _{n}\ll n^{\gamma _{i}}%
\big|\,\mathbf{Z}(n)\neq \mathbf{0}\right) =1-\frac{1}{2^{i-1}};
\label{recent_i}
\end{equation}%
3a) for $i=1,2,...,N-1$
\begin{equation}
\lim_{n\rightarrow \infty }\mathbf{P}\left( n^{\gamma _{i}}\ll \beta _{n}\ll
n^{\gamma _{i+1}}\big|\,\mathbf{Z}(n)\neq \mathbf{0}\right) =0;
\label{NoMRCA}
\end{equation}%
4) for any $x\in (0,1)$
\begin{equation*}
\lim_{n\rightarrow \infty }\mathbf{P}(\beta _{n}\leq xn|\mathbf{Z}(n)\neq
\mathbf{0})=1-\frac{1}{2^{N-1}}(1-x).
\end{equation*}
\end{theorem}

\textbf{Remark 3.} As we see by (\ref{NoMRCA}), there are time-intervals of
increasing orders within each of which the probability to find the MRCA of
the population survived up to moment $n\rightarrow \infty $ is negligible
compared to the probability for the population to survive up to this moment.
Moreover, these time-intervals are separated from each other by the
time-intervals of increasing orders within each of which the probability to
find the MRCA is strictly positive. Such a phenomena has no analogues for
the indecomposable Galton-Watson processes.

Along with the distribution of the birth moment of the MRCA, the type $\zeta
_{n}$ of the MRCA of the population survived up to moment $n$ is of
interest. The distribution of this random variable is described by the
following theorem.

\begin{theorem}
\label{T_type}Let Hypothesis A be valid. Then, for $i=1,2,...,N$
\begin{equation*}
p_{i}=\lim_{n\rightarrow \infty }\mathbf{P}(\zeta _{n}=i|\mathbf{Z}(n)\neq
\mathbf{0})=\frac{1}{2^{i}}(1-\delta _{iN})+\frac{1}{2^{N-1}}\,\delta _{iN},
\end{equation*}%
where $\delta _{ij}$ is the Kroneker symbol.
\end{theorem}

Observe that $p_{N-1}=p_{N}$.

\textbf{Remark 4.} The authors of paper \cite{FN2}, which contains several
results used in the proofs of our Theorems \ref{T_SkorohConst}-\ref{T_type},
considered a more general case of the strongly critical branching processes.
Namely, they prove a number of conditional limit theorems for the case when
by a suitable labelling the types of the multitype Galton-Watson process can
be grouped into $N\geq 2$ partially ordered classes $\mathcal{C}%
_{1}\rightarrow \mathcal{C}_{2}\rightarrow ...\rightarrow \mathcal{C}_{N}$
possessing the following properties:

1) particle types belonging to any given class, say $\mathcal{C}_{i},$
constitute an indecomposable critical branching process with $r_{i}\geq 1$
types;

2) each class $\mathcal{C}_{i}$ contains a type whose representatives are
able to produce offspring in the next class in the order with a positive
probability;

3) particles with types from $\mathcal{C}_{i},i\geq 2,$ are unable to
produce offspring belonging to the classes $\mathcal{C}_{1},...,\mathcal{C}%
_{i-1}$.

The methods used in the present paper may be applied to investigate, for
instance, the asymptotic distribution of $\beta _{n}$ for such processes.
Since the needed arguments are too cumbersome and contain no new ideas, we
prefer to concentrate on the case when each class $\mathcal{C}_{i}$ includes
a single type only.

The remainder of the paper is organized as follows. Section \ref{Sec2}
contains some preliminary results. In particular, we recall the statements
from \cite{FN} and \cite{FN2} describing the asymptotic behavior of the
survival probability and the distribution of the number of particles in a
strongly critical decomposable branching process. Section \ref{Sec3} gives a
detailed description of the limiting processes. In Sections \ref{Sec4} and %
\ref{Sec5} we check convergence of one-dimensional and finite-dimensional
distributions of the prelimiting processes to the limiting ones. Section \ref%
{Sec6} contains the proofs of Theorems \ref{T_SkorohConst} and \ref{T_Skhod1}%
. Finally, Section \ref{Sec7} is devoted to the proofs of Theorems \ref%
{T_mrcaMany} and \ref{T_type}.

\section{Auxiliary results\label{Sec2}}

For any vector $\mathbf{s}=(s_{1},...,s_{p})$ (the dimension will usually be
clear from the context) and an integer valued vector $\mathbf{k}%
=(k_{1}.....k_{p})$ define%
\begin{equation*}
\mathbf{s}^{\mathbf{k}}=s_{1}^{k_{1}},...,\,s_{p}^{k_{p}}.
\end{equation*}%
Further, let $\mathbf{1}=(1,...,1)$ be a vector of units. It will be
sometimes convenient to write $\mathbf{1}^{(i)}$ for the $i-$dimensional
vector with all its components equal to one.

Let
\begin{equation*}
H_{n}^{(i,N)}(\mathbf{s})=\mathbf{E}\left[ \mathbf{s}^{\mathbf{Z}(n)}|%
\mathbf{Z}(0)=\mathbf{e}_{i}\right] =\mathbf{E}\left[ s_{i}^{Z_{i}(n)}...%
\,s_{N}^{Z_{N}(n)}|\mathbf{Z}(0)=\mathbf{e}_{i}\right]
\end{equation*}%
be the probability generating function for $\mathbf{Z}(n)$ given the process
is initiated at time zero by a single particle of type $i\in \left\{
1,2,...,N\right\}.$ Clearly (recall (\ref{DefNONimmigr})), $H_{1}^{(i,N)}(%
\mathbf{s})=h_{i}(\mathbf{s}),\ i=1,...,N$. Denote
\begin{equation*}
Q_{n}^{(i,N)}(\mathbf{s})=1-H_{n}^{(i,N)}(\mathbf{s}%
),~Q_{n}^{(i,N)}=1-H_{n}^{(i,N)}(\mathbf{0}),
\end{equation*}
put
\begin{equation*}
\mathbf{H}_{n}(\mathbf{s})=(H_{n}^{(1,N)}(\mathbf{s}),...,H_{n}^{(N,N)}(%
\mathbf{s})),~\mathbf{Q}_{n}(\mathbf{s})=(Q_{n}^{(1,N)}(\mathbf{s}%
),...,Q_{n}^{(N,N)}(\mathbf{s}))
\end{equation*}%
and set%
\begin{equation*}
b_{jk}(n)=\mathbf{E}\left[ Z_{j}(n)Z_{k}(n)-\delta _{jk}Z_{j}(n)|\mathbf{Z}%
(0)=\vec{e}_{j}\right].
\end{equation*}

The starting point of our arguments is the following theorem being a
simplified combination of the respective results from \cite{FN} and \cite%
{FN2}:

\begin{theorem}
\label{T_Foster}Let $\mathbf{Z}(n),n=0,1,..$ be a strongly critical
decomposable multitype branching process satisfying (\ref{matrix1}), (\ref%
{Matpos}), (\ref{Maseq}), and (\ref{FinCovar}). Then, as $n\rightarrow
\infty $%
\begin{eqnarray}
m_{jj}(n) &=&1,\ m_{ij}(n)\sim a_{ij}n^{j-i},\ i<j,  \label{MomentSingle3} \\
b_{jk}(n) &\sim &\hat{a}_{jk}n^{k-j+1},\ j\leq k,  \label{Momvariance}
\end{eqnarray}%
where $a_{ij}$ and $\hat{a}_{jk}$ are positive constants known explicitly
(see \cite{FN2}, Theorem 1).

Besides \textit{(}see \cite{FN}, Theorem 1), as $n\rightarrow \infty $
\begin{equation}
Q_{n}^{(i,N)}=1-H_{n}^{(i,N)}(\mathbf{0})=\mathbf{P}(\mathbf{Z}( n) \neq
\mathbf{0}|\mathbf{Z}( 0) =\mathbf{e}_{i})\sim c_{iN}n^{-1/2^{N-i}},
\label{SurvivSingle}
\end{equation}%
where the constants $c_{iN}$ are the same as in (\ref{Const2}).
\end{theorem}

In the sequel we prove the following Yaglom-type limit theorem being a
compliment to Theorem \ref{T_Foster}.

\begin{theorem}
\label{T_Yaglom}Under the conditions of Theorem \ref{T_Foster}, for any $%
\lambda >0$%
\begin{equation}
\lim_{n\rightarrow \infty }\mathbf{E}\left[ \exp \left\{ -\lambda \frac{%
Z_{N}(n)}{b_{N}n}\right\} \Big|\mathbf{Z}( n) \neq \mathbf{0};\mathbf{Z}( 0)
=\mathbf{e}_{i}\right] =1-\Big( \frac{\lambda }{1+\lambda }\Big) %
^{1/2^{N-i}}.  \label{Yag}
\end{equation}
\end{theorem}

Set $d_{ii}=\sqrt{b_{i}^{-1}m_{i,i+1}}$ , $i=1,2,...,N-1$\ and, for $%
j=1,2,...,i-1$ let%
\begin{equation}
d_{ji}=\sqrt{b_{j}^{-1}m_{j,j+1}d_{j+1,i}},\quad D_{i}=d_{1i}.
\label{Dreccur}
\end{equation}%
Observe that (see (\ref{Const1})) for $k=0,1,2,...,i-1$
\begin{equation}
d_{i-k,i}=(b_{i}m_{i,i+1})^{1/2^{k+1}}c_{i-k,i},\quad
D_{i}=(b_{i}m_{i,i+1})^{1/2^{i}}c_{1i}=(b_{i}m_{i,i+1})^{1/2^{i}}C_{i}.
\label{DcConnection}
\end{equation}

Let $\mathbf{Z}(0)=\mathbf{e_{1}}$ and denote by
\begin{equation*}
T_{i}=\min \left\{ n\geq 1:Z_{1}(n)+Z_{2}(n)+...+Z_{i}(n)=0\right\}
\end{equation*}%
the extinction moment of the population generated by the particles of the
first $i$ in order types. Let $\eta _{rj}\left( k,l\right) $ be the number
of daughters of type $j$ of the $l-$th mother of type $r$ belonging to the $%
k-$th generation and
\begin{equation*}
W_{pij}=\sum_{r=p}^{i}\sum_{k=0}^{T_{i}}\sum_{q=1}^{Z_{r}(k)}\eta
_{rj}\left( k,q\right)
\end{equation*}%
be the total amount of daughters of type $j\geq i+1$ produced by all
particles of types $p,p+1,...,i$ ever born in the process if the process is
initiated at time $n=0$ by a single particle of type $p\leq i.$ Finally, put%
\begin{equation*}
W_{pi}=\sum_{j=i+1}^{N}W_{pij}=\sum_{j=i+1}^{N}\sum_{r=p}^{i}%
\sum_{k=0}^{T_{i}}\sum_{q=1}^{Z_{r}(k)}\eta _{rj}\left( k,q\right) .
\end{equation*}

We know by (\ref{SurvivSingle}) that%
\begin{equation}
Q_{n}^{(1,i)}=\mathbf{P}( T_{i}>n) \sim c_{1i}n^{-2^{-(i-1)}}.  \label{Asqq}
\end{equation}

The next lemma describes the tail distributions of $W_{1i,i+1}$ and $W_{1i}$.

\begin{lemma}
\label{L_Laplace} Let Hypothesis A be valid. Then, as $\lambda \downarrow 0$%
\begin{equation}
1-\mathbf{E}\left[ e^{-\lambda W_{1i,i+1}}\,|\mathbf{Z}(0)=\mathbf{e}_{1}%
\right] \sim d_{1i}\lambda ^{1/2^{i}}=D_{i}\lambda ^{1/2^{i}}  \label{Tot1}
\end{equation}%
and there exists a constant $F_{i}>0$ such that%
\begin{equation}
1-\mathbf{E}\left[ e^{-\lambda W_{1i}}|\mathbf{Z}(0)=\mathbf{e}_{1}\right]
\sim F_{i}\lambda ^{1/2^{i}}.  \label{Tot2}
\end{equation}
\end{lemma}

\textbf{Proof.} Set
\begin{equation*}
W_{pi,i+1}(n)=\sum_{r=p}^{i}\sum_{k=0}^{n}\sum_{q=1}^{Z_{r}(k)}\eta
_{rj}\left( k,q\right) ,
\end{equation*}%
denote%
\begin{equation*}
K_{pi,n}(\mathbf{s};t)=\mathbf{E}\left[
s_{p}^{Z_{p}(n)}...s_{i}^{Z_{i}(n)}t^{W_{pi,i+1}(n)}|\mathbf{Z}(0)=\mathbf{e}%
_{p}\right] ,\,K_{pi,n}(t)=K_{pi,n}(\mathbf{1}^{(i-p+1)};t)
\end{equation*}%
and put
\begin{equation*}
K_{pi}(t)=\mathbf{E}\left[ t^{W_{pi,i+1}}\big|\,\mathbf{Z}(0)=\mathbf{e}_{p}%
\right] =\lim_{n\rightarrow \infty }K_{pi,n}(t)
\end{equation*}%
(this limit exists since the random variables $W_{pi,i+1}(n),p=1,2,...,i$
are nondecreasing in $n$). Clearly, to prove the lemma it is sufficient to
show that, as~$t\uparrow 1$%
\begin{equation*}
1-K_{1i}(t)=1-\mathbf{E}\left[ t^{W_{1i,i+1}}\,|\mathbf{Z}(0)=\mathbf{e}_{1}%
\right] {}\mathbf{\sim {}}d_{1i}(1-t)^{1/2^{i}}.
\end{equation*}

Using properties of branching processes it is not difficult to check that%
\begin{equation*}
K_{pi,n+1}(\mathbf{s};t)=h_{p}\left( K_{pi,n}(\mathbf{s};t),...,K_{ii,n}(%
\mathbf{s};t),t,\mathbf{1}^{(N-i-1)}\right)
\end{equation*}%
implying%
\begin{equation*}
K_{pi,n+1}(t)=h_{p}\left( K_{pi,n}(t),...,K_{ii,n}(t),t,\mathbf{1}%
^{(N-i-1)}\right) .
\end{equation*}%
and%
\begin{equation*}
K_{pi}(t)=h_{p}\left( K_{pi}(t),...,K_{ii}(t),t,\mathbf{1}^{(N-i-1)}\right) .
\end{equation*}%
In particular,%
\begin{equation*}
K_{ii}(t)=h_{i}\left( K_{ii}(t),t,\mathbf{1}^{(N-i-1)}\right) .
\end{equation*}%
Since $\mathbf{E}\eta _{ii}=1$ and $b_{i}=\frac{1}{2}Var\eta _{ii}\in
(0,\infty )$, it follows that, as $t\uparrow 1$%
\begin{eqnarray*}
1-K_{ii}(t) &=&1-h_{i}\left( K_{ii}(t),t,\mathbf{1}^{(N-i-1)}\right) \\
&=&1-K_{ii}(t)-b_{i}(1-K_{ii}(t))^{2}(1+o(1))+m_{i,i+1}(1-t)
\end{eqnarray*}%
or%
\begin{equation*}
1-K_{ii}(t)\sim \sqrt{b_{i}^{-1}m_{i,i+1}(1-t)}.
\end{equation*}%
This, in particular, proves the statement of the lemma for $i=1$.

Now we use induction and assume that%
\begin{equation*}
1-K_{qi}(t)\sim d_{qi}(1-t)^{1/2^{i-q+1}},q=p+1,...,i.
\end{equation*}%
Then%
\begin{eqnarray*}
1-K_{pi}(t) &=&1-h_{p}\left( K_{pi}(t),...,K_{ii}(t),t,\mathbf{1}%
^{(N-i-1)}\right) \\
&=&1-K_{pi}(t)-b_{p}(1-K_{pi}(t))^{2}(1+o(1)) \\
&&+(1+o(1))\left( m_{p,p+1}(1-K_{p+1,i}(t))+\sum_{q=p+2}^{i}m_{pq}\left(
1-K_{qi}(t)\right) \right) \\
&&+(1+o(1))m_{p,i+1}(1-t)
\end{eqnarray*}%
implying%
\begin{eqnarray*}
1-K_{pi}(t) &\sim &\sqrt{b_{p}^{-1}m_{p,p+1}(1-K_{p+1,i}(t))} \\
&\sim &\sqrt{b_{p}^{-1}m_{p,p+1}d_{p+1,i}}\left( 1-t\right)
^{1/2^{i-p+1}}=d_{pi}\left( 1-t\right) ^{1/2^{i-p+1}}
\end{eqnarray*}%
and proving (\ref{Tot1}).

To prove (\ref{Tot2}) it is necessary to use similar arguments. We omit the
details.

Lemma \ref{L_Laplace} is proved.\

From now on and till the end of this section we suppose that%
\begin{equation}
s_{k}=\exp (-\lambda _{k}n^{-2^{-(N-k)}})=\exp (-\lambda _{k}n^{-\gamma
_{k}}),\lambda _{k}>0,\,k=1,2,...,N  \label{Sasymp}
\end{equation}%
and, keeping in mind this assumption, study in Lemmas \ref{L_OneOnly}-\ref%
{L_MultiSharp} the asymptotic behavior of the difference $1-H_{m}^{(j,N)}(%
\mathbf{s})$ when $m,n\rightarrow \infty .$

\begin{lemma}
\label{L_OneOnly}If%
\begin{equation}
m\ll n^{2^{-(N-j)}}=n^{\gamma _{j}}  \label{Mnegl}
\end{equation}%
then for $N>j$%
\begin{equation*}
\lim_{n\rightarrow \infty }n^{\gamma _{j}}Q_{m}^{(j,N)}(\mathbf{s})=\lambda
_{j}.
\end{equation*}
\end{lemma}

\textbf{Proof.} Clearly, it is sufficient to prove the statement for $j=1$
only. Let $r$ be a positive integer such that%
\begin{equation*}
1-H_{r}^{(1,1)}(0)\leq 1-s_{1}\leq 1-H_{r-1}^{(1,1)}(0).
\end{equation*}%
Since $1-s_{1}\sim \lambda _{1}n^{-\gamma _{1}}$ and $1-H_{r}^{(1,1)}(0)\sim
(b_{1}r)^{-1}\,$as $n,r\rightarrow \infty $, it follows that $r\sim
(b_{1}\lambda _{1})^{-1}n^{\gamma _{1}}$. By the branching property of
probability generating functions we have for $m\ll n^{\gamma _{1}}:$
\begin{eqnarray*}
Q_{m}^{(1,N)}(\mathbf{s}) &\geq &1-H_{m}^{(1,1)}(s_{1})\geq
1-H_{m}^{(1,1)}(H_{r}^{(1,1)}(0)) \\
&=&1-H_{m+r}^{(1,1)}(0)\sim b_{1}^{-1}(m+r)^{-1}\sim \lambda _{1}n^{-\gamma
_{1}}.
\end{eqnarray*}%
Besides,%
\begin{eqnarray*}
Q_{m}^{(1,N)}(\mathbf{s}) &\leq &1-H_{m}^{(1,1)}(s_{1})+\mathbf{E}\left[
\left( 1-s_{2}^{Z_{2}(m)}...\,s_{N}^{Z_{N}(m)}\right) |\mathbf{Z}(0)=\mathbf{%
e}_{1}\right] \\
&\leq &1-H_{m+r-1}^{(1,1)}(0)+\sum_{k=2}^{N}(1-s_{k})\mathbf{E}\left[
Z_{k}(m)|\mathbf{Z}(0)=\mathbf{e}_{1}\right] .
\end{eqnarray*}%
We know by (\ref{MomentSingle3}) and (\ref{Sasymp}) that, for a positive
constant $C$
\begin{equation*}
\sum_{k=2}^{N}(1-s_{k})\mathbf{E}\left[ Z_{k}(m)|\mathbf{Z}(0)=\mathbf{e}_{1}%
\right] \leq C\sum_{k=2}^{N}\lambda _{k}n^{-\gamma _{k}}m^{k-1}
\end{equation*}%
which, in view of (\ref{Mnegl}) is negligible with respect to%
\begin{equation*}
C\max_{2\leq i\leq N}\lambda _{i}\times \sum_{k=2}^{N}n^{-\gamma
_{k}}(n^{\gamma _{1}})^{k-1}=C\max_{2\leq i\leq N}\lambda _{i}\times
\sum_{k=2}^{N}n^{(k-1)2^{-(N-1)}-2^{-(N-k)}}.
\end{equation*}%
Since $k2^{-(N-1)}-2^{-(N-k)}=2^{-(N-1)}(k-2^{k-1})\leq 0$ for $k\geq 2,$ we
have
\begin{equation*}
n^{2^{-(N-1)}}\sum_{k=2}^{N}n^{(k-1)2^{-(N-1)}-2^{-(N-k)}}=%
\sum_{k=2}^{N}n^{k2^{-(N-1)}-2^{-(N-k)}}\leq N-1.
\end{equation*}%
Consequently, $Q_{m}^{(1,N)}(\mathbf{s})\sim 1-H_{m}^{(1,1)}(s_{1})\sim
\lambda _{1}n^{-\gamma _{1}}$ as $n\rightarrow \infty $.

This proves the lemma.

In order to formulate the next lemma we introduce a tuple of functions $\phi
_{i}=\phi _{i}( \lambda _{1},\lambda _{2}) ,\ i=1,2,...,N-1$ solving in the
domain $\left\{ \lambda _{1}\geq 0,\lambda _{2}\geq 0\right\} $ the
differential equations%
\begin{equation*}
\lambda _{1}\frac{\partial \phi _{i}}{\partial \lambda _{1}}+2\lambda _{2}%
\frac{\partial \phi _{i}}{\partial \lambda _{2}}=-b_{i}\phi _{i}^{2}+\phi
_{i}+m_{i,i+1}\lambda _{2}
\end{equation*}%
with the initial conditions
\begin{equation*}
\phi _{i}(\mathbf{0})=0,\ \frac{\partial \phi _{i}( \mathbf{0}) }{\partial
\lambda _{1}}=1,\ \frac{\partial \phi _{i}( \mathbf{0}) }{\partial \lambda
_{2}}=m_{i,i+1}.
\end{equation*}

One may check that, for any $y>0$
\begin{equation}
\frac{\phi _{i}(\lambda _{1}y,\lambda _{2}y^{2})}{y}=\sqrt{\frac{%
m_{i,i+1}\lambda _{2}}{b_{i}}}\frac{b_{i}\lambda _{1}+\sqrt{%
b_{i}m_{i,i+1}\lambda _{2}}\tanh (y\sqrt{b_{i}m_{i,i+1}\lambda _{2}})}{%
b_{i}\lambda _{1}\tanh (y\sqrt{b_{i}m_{i,i+1}\lambda _{2}})+\sqrt{%
b_{i}m_{i,i+1}\lambda _{2}}}.  \label{DefSimpl}
\end{equation}

\begin{lemma}
\label{L_twoOnly}Let condition (\ref{Sasymp}) be valid. If $m\sim yn^{\gamma
_{i}},$ $y>0$ then%
\begin{equation*}
\lim_{n\rightarrow \infty }n^{\gamma _{i}}Q_{m}^{(i,N)}(\mathbf{s}%
)=y^{-1}\phi _{i}(\lambda _{i}y,\lambda _{i+1}y^{2}).
\end{equation*}
\end{lemma}

\textbf{Proof.} As in the previous lemma, it is sufficient to consider the
case $i=1$ only. It follows from Theorem 2 in \cite{FN2} that for $\lambda
_{k}\geq 0,k=1,2,...,N$%
\begin{eqnarray*}
&&\lim_{m\rightarrow \infty }m\left( 1-\mathbf{E}\left[ \exp \left\{
-\sum_{k=1}^{N}\lambda _{k}\frac{Z_{k}(m)}{m^{k}}\right\} \right] \right) \\
&&\qquad \qquad =\lim_{m\rightarrow \infty }m(1-H_{m}^{(1,N)}(e^{-\lambda
_{1}/m},e^{-\lambda _{2}/m^{2}},...,e^{-\lambda _{N}/m^{N}})) \\
&&\qquad \qquad =\Phi (\lambda _{1},\lambda _{2},...,\lambda _{N}),
\end{eqnarray*}%
where $\Phi =\Phi (\lambda _{1},\lambda _{2},...,\lambda _{N})$ solves the
differential equation%
\begin{equation*}
\sum_{k=1}^{N}k\lambda _{k}\frac{\partial \Phi }{\partial \lambda _{k}}%
=-b_{1}\Phi ^{2}+\Phi +\sum_{k=2}^{N}f_{k}\lambda _{k}
\end{equation*}%
with the initial conditions%
\begin{equation*}
\Phi (\mathbf{0})=0,\ \frac{\partial \Phi (\mathbf{0})}{\partial \lambda _{1}%
}=1,\ \frac{\partial \Phi (\mathbf{0})}{\partial \lambda _{k}}=\frac{1}{k-1}%
f_{k},\ k=2,...,N
\end{equation*}%
and
\begin{equation*}
f_{k}=\frac{1}{(k-2)!}\prod_{j=1}^{k-1}m_{j,j+1},\ k=2,...,N.
\end{equation*}%
Since $m^{2^{k-1}}=m^{k}$ for $k=1,2$ and $m^{2^{k-1}}\gg m^{k}$ for $k>2,$
we conclude by the continuity of $\Phi $ at point $\mathbf{0}$ that%
\begin{eqnarray*}
&&\lim_{n\rightarrow \infty }n^{\gamma _{1}}Q_{m}^{(1,N)}(\mathbf{s}%
)=y^{-1}\lim_{m\rightarrow \infty }mQ_{m}^{(1,N)}(\mathbf{s}) \\
&&\quad =y^{-1}\lim_{m\rightarrow \infty }m\left( 1-\mathbf{E}\left[ \exp
\left\{ -\sum_{k=1}^{N}\lambda _{k}\frac{Z_{k}(m)}{n^{1/2^{N-k}}}\right\} %
\right] \right) \\
&&\quad =y^{-1}\lim_{m\rightarrow \infty }m\left( 1-\mathbf{E}\left[ \exp
\left\{ -\sum_{k=1}^{N}\lambda _{k}y^{2^{k-1}}\frac{Z_{1}(m)}{m^{2^{k-1}}}%
\right\} \right] \right) \\
&&\quad =y^{-1}\Phi (\lambda _{1}y,\lambda _{2}y^{2},0,...,0)=y^{-1}\phi
_{1}(\lambda _{1}y,\lambda _{2}y^{2}).
\end{eqnarray*}

Lemma \ref{L_twoOnly} is proved.\

\begin{lemma}
\label{L_inbetween}Let condition (\ref{Sasymp}) be valid. If, for some $%
i\leq N-1$
\begin{equation}
n^{\gamma _{i}}\ll m\ll n^{\gamma _{i+1}}  \label{msmall}
\end{equation}%
then
\begin{equation*}
\lim_{n\rightarrow \infty }n^{\gamma _{1}}Q_{m}^{(1,N)}(\mathbf{s}%
)=D_{i}\left( \lambda _{i+1}\right) ^{1/2^{i}}.
\end{equation*}
\end{lemma}

\textbf{Proof.} It follows from (\ref{Asqq}) and (\ref{msmall}) that
\begin{equation*}
\mathbf{P}(T_{i}>m)\sim c_{1i}m^{-2^{-(i-1)}}=o(n^{-\gamma _{1}}).
\end{equation*}%
Therefore,%
\begin{eqnarray*}
Q_{m}^{(1,N)}(\mathbf{s}) &=&\mathbf{E}\left[
1-s_{1}^{Z_{1}(m)}s_{2}^{Z_{2}(m)}...\,s_{N}^{Z_{N}(m)}\right] \\
&=&\mathbf{E}\left[ \left(
1-s_{i+1}^{Z_{i+1}(m)}...\,s_{N}^{Z_{N}(m)}\right) ;T_{i}\leq m\right]
+o(n^{-\gamma _{1}}) \\
&=&1-H_{m}^{(1,N)}\left( \mathbf{1}^{(i)},s_{i+1},...,s_{N}\right)
+o(n^{-\gamma _{1}}).
\end{eqnarray*}%
It is not difficult to check that for our decomposable branching process
\begin{eqnarray*}
&&H_{m}^{(1,N)}\left( \mathbf{1}^{(i)},s_{i+1},...,s_{N}\right) \\
&&\quad =\mathbf{E}\left[ \prod_{k=0}^{m-1}\prod_{r=1}^{i}%
\prod_{l=1}^{Z_{r}(k)}\prod_{j=i+1}^{N}\left( H_{m-k}^{(j,N)}(\mathbf{s}%
)\right) ^{\eta _{rj}\left( k,l\right) }\right] \\
&&\quad =\mathbf{E}\left[ \prod_{k=0}^{m-1}\prod_{r=1}^{i}%
\prod_{l=1}^{Z_{r}(k)}\prod_{j=i+1}^{N}\left( H_{m-k}^{(j,N)}(\mathbf{s}%
)\right) ^{\eta _{rj}\left( k,l\right) };T_{i}\leq \sqrt{mn^{\gamma _{i}}}%
\right] \\
&&\qquad +O\left( \mathbf{P}\left( T_{i}> \sqrt{mn^{\gamma _{i}}}\right)
\right) .
\end{eqnarray*}

Observing that $\lim_{m\to\infty}H_{m-k}^{(j,N)}(\mathbf{s})\rightarrow 1$
for $j\geq i+1$ and $k\leq T_{i}\leq \sqrt{mn^{\gamma _{i}}}=o(m),$ we get
on the set $T_{i}\leq \sqrt{mn^{\gamma _{i}}}$
\begin{eqnarray*}
&&\prod_{k=0}^{m-1}\prod_{r=1}^{i}\prod_{l=1}^{Z_{r}(k)}\prod_{j=i+1}^{N}%
\left( H_{m-k}^{(j,N)}(\mathbf{s})\right) ^{\eta _{rj}\left( k,l\right) } \\
&&\quad =\exp \left\{
-\sum_{r=1}^{i}\sum_{k=0}^{T_{i}}\sum_{l=1}^{Z_{r}(k)}\sum_{j=i+1}^{N}\eta
_{rj}\left( k,l\right) Q_{m-k}^{(j,N)}(\mathbf{s})(1+o(1))\right\} .
\end{eqnarray*}%
If $j\geq i+1$ then Lemma \ref{L_OneOnly} and the estimates $m\ll n^{\gamma
_{i+1}}\leq n^{\gamma _{j}}$ yield
\begin{equation*}
Q_{m-k}^{(j,N)}(\mathbf{s})\sim Q_{m}^{(j,N)}(\mathbf{s})\sim \lambda
_{j}n^{-\gamma _{j}}.
\end{equation*}%
Hence it follows that on the set $T_{i}\leq \sqrt{mn^{\gamma _{i}}}%
=o(m)=o(n^{\gamma _{i+1}})$%
\begin{eqnarray*}
&&\sum_{r=1}^{i}\sum_{k=0}^{T_{i}}\sum_{l=1}^{Z_{r}(k)}\sum_{j=i+1}^{N}\eta
_{rj}\left( k,l\right) Q_{m-k}^{(j,N)}(\mathbf{s}) \\
&&\quad =(1+o(1))\sum_{j=i+1}^{N}Q_{m}^{(j,N)}(\mathbf{s})\sum_{r=1}^{i}%
\sum_{k=0}^{T_{i}}\sum_{l=1}^{Z_{r}(k)}\eta _{rj}\left( k,l\right) \\
&&\quad =(1+o(1))\sum_{j=i+1}^{N}W_{1ij}Q_{m}^{(j,N)}(\mathbf{s}) \\
&&\quad =(1+o(1))W_{1i,i+1}Q_{m}^{(i+1,N)}(\mathbf{s})+O\left(
Q_{m}^{(i+2,N)}(\mathbf{s})\right) \sum_{j=i+2}^{N}W_{1ij} \\
&&\quad =(1+o(1))W_{1i,i+1}\lambda _{i+1}n^{-\gamma _{i+1}}+O_{n}(n^{-\gamma
_{i+2}}W_{1i}).
\end{eqnarray*}%
Using the estimates
\begin{eqnarray*}
0 &\leq &\mathbf{E}\left[ \exp \left\{ -(1+o(1))W_{1i,i+1}\lambda
_{i+1}n^{-\gamma _{i+1}}\right\} \right] \\
&&-\mathbf{E}\left[ \exp \left\{ -(1+o(1))W_{1i,i+1}\lambda _{i+1}n^{-\gamma
_{i+1}}-O(n^{-\gamma _{i+2}}W_{1i})\right\} \right] \\
&\leq &1-\mathbf{E}\left[ \exp \left\{ -O(n^{-\gamma _{i+2}}W_{1i})\right\} %
\right] =O\left( \left( n^{-\gamma _{i+2}}\right) ^{1/2^{i}}\right) =O\left(
n^{-\gamma _{2}}\right)
\end{eqnarray*}%
where, for the penultimate equality we applied (\ref{Tot2}), we conclude by (%
\ref{Tot1}) that%
\begin{eqnarray*}
&&1-H_{m}^{(1,N)}\left( \mathbf{1}^{(i)},s_{i+1},...,s_{N}\right) \\
&&\quad =(1+o(1))\mathbf{E}\left[ 1-\exp \left\{ -(1+o(1))W_{1i,i+1}\lambda
_{i+1}n^{-\gamma _{i+1}}\right\} \right] \\
&&\qquad +O\left( \mathbf{P}\left( T_{i}>\sqrt{mn^{\gamma _{i}}}\right)
\right) \\
&&\quad =(1+o(1))D_{i}\left( \lambda _{i+1}n^{-\gamma _{i+1}}\right)
^{1/2^{i}}+o(n^{-\gamma _{1}})\sim D_{i}(\lambda _{i+1})^{1/2^{i}}n^{-\gamma
_{1}}
\end{eqnarray*}%
as desired. \

\begin{lemma}
\label{L_MultiSharp}If $m\sim yn^{\gamma _{i}}$ for some $i\in \left\{
2,3,...,N-1\right\} $ then%
\begin{equation*}
\lim_{n\rightarrow \infty }n^{\gamma _{1}}Q_{m}^{(1,N)}(\mathbf{s}%
)=D_{i-1}(y^{-1}\phi _{i}(\lambda _{i}y,\lambda _{i+1}y^{2}))^{1/2^{i-1}}.
\end{equation*}
\end{lemma}

\textbf{Proof.} If $m\sim yn^{\gamma _{i}}$ and $j\geq i$ then $n^{\gamma
_{j}}\sim (y^{-1}m)^{2^{j-i}}$ and, therefore,%
\begin{equation*}
s_{j}=\exp \left\{ -\lambda _{j}n^{-\gamma _{j}}\right\} =\exp \left\{
-(1+o(1))\lambda _{j}y^{2^{j-i}}m^{-2^{j-i}}\right\} .
\end{equation*}%
Hence we may apply Lemma \ref{L_twoOnly} to get, as $n\rightarrow \infty $%
\begin{equation*}
n^{\gamma _{i}}Q_{m}^{(i,N)}(\mathbf{s})\sim
y^{-1}mQ_{m}^{(i,N)}(s_{i},s_{i+1},...,s_{N})\sim y^{-1}\phi _{i}(\lambda
_{i}y,\lambda _{i+1}y^{2}).
\end{equation*}%
Further, as in the previous lemma we have
\begin{equation*}
Q_{m}^{(1,N)}(\mathbf{s})=1-H_{m}^{(1,N)}\left( \mathbf{1}%
^{(i-1)},s_{i},...,s_{N}\right) +o(n^{-\gamma _{1}})
\end{equation*}%
and on the set $T_{i-1}\leq \sqrt{mn^{\gamma _{i-1}}}\ll m\sim yn^{\gamma
_{i}}$%
\begin{eqnarray*}
&&\sum_{r=1}^{i-1}\sum_{k=0}^{T_{i-1}}\sum_{l=1}^{Z_{r}(k)}\sum_{j=i}^{N}%
\eta _{rj}\left( k,l\right) Q_{m-k}^{(j,N)}(\mathbf{s}) \\
&&\quad \quad =(1+o(1))\sum_{j=i}^{N}W_{1,i-1,j}Q_{m}^{(j,N)}(\mathbf{s}) \\
&&\quad =(1+o(1))W_{1,i-1,i}Q_{m}^{(i,N)}(\mathbf{s})+O\left(
Q_{m}^{(i+1,N)}(\mathbf{s})\right) \sum_{j=i+1}^{N}W_{1,i-1,j} \\
&&\quad =(1+o(1))W_{1,i-1,i}(y^{-1}\phi _{i}(\lambda _{i}y,\lambda
_{i+1}y^{2}))^{1/2^{i-1}}n^{-\gamma _{i+1}} \\
&&+O_{n}(n^{-\gamma _{i+2}}W_{1,i-1}).
\end{eqnarray*}%
Therefore,%
\begin{eqnarray*}
&&1-H^{(1,N)}_{m}\left( \mathbf{1}^{(i-1)},s_{i},...,s_{N}\right) \\
&&\quad =\mathbf{E}\left[ 1-\exp \left\{ -(1+o(1))W_{1,i-1,i}y^{-1}\phi
_{i}(\lambda _{i}y,\lambda _{i+1}y^{2})\,n^{-\gamma _{i}}\right\} \right] \\
&&\qquad +O\left( \mathbf{P}\left( T_{i-1}\geq \sqrt{mn^{\gamma _{i+1}}}%
\right) \right) \\
&&\quad =(1+o(1))D_{i-1}\big(y^{-1}\phi _{i}(\lambda _{i}y,\lambda
_{i+1}y^{2})\,n^{-\gamma _{i}}\big)^{1/2^{i-1}}+o(n^{-\gamma _{1}}) \\
&&\quad \sim D_{i-1}(y^{-1}\phi _{i}(\lambda _{i}y,\lambda
_{i+1}y^{2}))^{1/2^{i-1}}n^{-\gamma _{1}}.
\end{eqnarray*}

The lemma is proved.

\begin{lemma}
\label{L_DC}For all $i=1,2,...,N-1$%
\begin{equation}
C_{N}=C_{i}(m_{i,i+1}b_{i}c_{i+1,N})^{1/2^{i}}=D_{i}(c_{i+1,N})^{1/2^{i}}.
\label{CD}
\end{equation}
\end{lemma}

\textbf{Proof.} Using (\ref{Const1}) we have
\begin{equation*}
c_{iN}=\sqrt{b_{i}^{-1}m_{i,i+1}c_{i+1,N}}=b_{i}^{-1}\sqrt{%
b_{i}m_{i,i+1}c_{i+1,N}}=c_{ii}\sqrt{b_{i}m_{i,i+1}c_{i+1,N}}
\end{equation*}%
leading in view of (\ref{Const2}) and (\ref{DcConnection}) to%
\begin{eqnarray*}
C_{N} &=&c_{1N}=\left( \frac{1}{b_{N}}\right)
^{1/2^{N-1}}\prod_{j=1}^{N-1}\left( \frac{m_{j,j+1}}{b_{j}}\right)
^{1/2^{j}}= \\
&=&c_{1i}(b_{i}m_{i,i+1})^{1/2^{i}}\left( \left( \frac{1}{b_{N}}\right)
^{1/2^{N-i}}\prod_{j=i+1}^{N-1}\left( \frac{m_{j,j+1}}{b_{j}}\right)
^{1/2^{j-i}}\right) ^{1/2^{i}} \\
&=&c_{1i}(b_{i}m_{i,i+1}c_{i+1,N})^{1/2^{i}}=D_{i}(c_{i+1,N})^{1/2^{i}}
\end{eqnarray*}%
as desired.

\section{Properties of the limiting processes\label{Sec3}}

In this section we give a more detailed description of the properties of the
limiting processes. It follows from the definition of $\mathbf{R}(t)$ that
if
\begin{equation*}
\mathbf{S}_{i}=(s_{i1},s_{i2},...,s_{iN})\in \left[ 0,1\right] ^{N}\text{
and }t_{i}\in \lbrack \gamma _{i-1},\gamma _{i}),i=1,2,...,N,
\end{equation*}%
then%
\begin{equation*}
\mathbf{E}\left[ \prod_{i=1}^{N}\mathbf{S}_{i}^{\mathbf{R}(t_{i})}\right]
=\Omega _{N}(s_{11},s_{22},...,s_{NN}),
\end{equation*}%
where $\Omega _{1}(s)=s$ and
\begin{equation}
\Omega _{i+1}(s_{1},s_{2},...,s_{i+1})=s_{1}\left( 1-\sqrt{1-\Omega
_{i}(s_{2},...,s_{i+1})}\right) ,\,i=1,2,....  \label{DefOmega}
\end{equation}

If now some intervals $[\gamma _{i-1},\gamma _{i})$ contain more than one
point of observation over the process $\mathbf{R}(\cdot )$, say, $\gamma
_{i-1}\leq t_{i1}<t_{i2}<...<t_{ik_{i}}<\gamma _{i},i=1,2,...,N,$ and $_{j}%
\mathbf{S}_{i}=\left( _{j}s_{i1},_{j}s_{i2},...,_{j}s_{iN}\right) \in \left[
0,1\right] ^{N}\text{ }$ then, clearly,%
\begin{equation*}
\mathbf{E}\left[ \prod_{i=1}^{N}\prod_{j=1}^{k_{i}}(_{j}\mathbf{S}_{i})^{%
\mathbf{R}(t_{ij})}\right] =\Omega _{N}\left(
\prod_{j=1}^{k_{1}}\,_{j}s_{11},\prod_{j=1}^{k_{2}}\,_{j}s_{22},...,%
\prod_{j=1}^{k_{N}}\,_{j}s_{NN}\right) .
\end{equation*}

To describe the characteristics of the processes $\mathbf{U}_{i}(\cdot
),i=1,...,N-1$, let, for $(s_{i},s_{i+1})\in \lbrack 0,1]^{2}$%
\begin{equation}
\varphi _{i}(y;s_{i},s_{i+1})=\sqrt{1-s_{i+1}}\frac{(1-s_{i})+\sqrt{1-s_{i+1}%
}\tanh (b_{i}c_{iN}y\sqrt{1-s_{i+1}})}{(1-s_{i})\tanh (b_{i}c_{iN}y\sqrt{%
1-s_{i+1}})+\sqrt{1-s_{i+1}}}  \label{Deffi}
\end{equation}%
with the natural agreement $\varphi _{i}(y;1,1)=0$ and%
\begin{equation*}
\varphi _{i}(y;s_{i},1)=\frac{1-s_{i}}{b_{i}c_{iN}y(1-s_{i})+1}.
\end{equation*}%
Denote%
\begin{eqnarray*}
X_{i}(y;\mathbf{s}) &=&X_{i}(y;s_{i},s_{i+1})=\mathbf{E}\left[ \mathbf{s}^{%
\mathbf{U}_{i}(y)}|\mathbf{U}_{i}(0)=\mathbf{e}_{i}\right] \\
&=&\mathbf{E}\left[ s_{i}^{U_{ii}(y)}s_{i+1}^{U_{i,i+1}(y)}|\mathbf{U}%
_{i}(0)=\mathbf{e}_{i}\right]
\end{eqnarray*}%
and set%
\begin{equation*}
\bar{X}_{R_{i}}(y;\mathbf{s})=\bar{X}_{R_{i}}(y;s_{i},s_{i+1})=\mathbf{E}%
_{R_{i}}\left[ \mathbf{s}^{\mathbf{U}_{i}(y)}\right] =\mathbf{E}_{R_{i}}%
\left[ s_{i}^{U_{ii}(y)}s_{i+1}^{U_{i,i+1}(y)}\right] ,
\end{equation*}%
where the symbol $\mathbf{E}_{R_{i}}[\cdot ]$ means that the process starts
by a random number of type $i$ particles distributed as $R_{i}$ in (\ref%
{Dist_ro}).

It follows from the description of the branching mechanism for $\mathbf{U}%
_{i}(\cdot )$ and the general theory of branching processes (see, for
instance, \cite{AN72}, p. 201) that $X_{i}(y;s_{i},s_{i+1})$ solves the
differential equation%
\begin{eqnarray*}
\frac{\partial }{\partial y}X_{i}(y;s_{i},s_{i+1}) &=&2b_{i}c_{iN}\left(
\frac{1}{2}X_{i}^{2}(y;s_{i},s_{i+1})-X_{i}(y;s_{i},s_{i+1})+\frac{1}{2}%
s_{i+1}\right) , \\
X_{i}(0;s_{i},s_{i+1}) &=&s_{i}.
\end{eqnarray*}%
Direct calculations show that%
\begin{equation}
X_{i}(y;s_{i},s_{i+1})=1-\varphi _{i}(y;s_{i},s_{i+1})  \label{ExplX}
\end{equation}%
and, as a result
\begin{equation}
\bar{X}_{R_{i}}(y;s_{i},s_{i+1})=1-(\varphi
_{i}(y;s_{i},s_{i+1}))^{1/2^{i-1}}\text{.}  \label{ExplXstar}
\end{equation}%
One may check by (\ref{Deffi}) and (\ref{ExplXstar}) that%
\begin{equation}
\lim_{y\downarrow 0}\bar{X}_{R_{i}}(y;s_{i},s_{i+1})=1-(1-s_{i})^{1/2^{i-1}}%
\text{ }  \label{LxLeft}
\end{equation}%
and%
\begin{equation}
\lim_{y\uparrow \infty }\bar{X}%
_{R_{i}}(y;s_{i},s_{i+1})=1-(1-s_{i+1})^{1/2^{i}}.  \label{liminf}
\end{equation}

For $y_{k}\in \lbrack 0,\infty ),\,(s_{ki},s_{k,i+1})\in \left[ 0,1\right]
^{2},\ k=1,2,...,p;\ i=1,...,N-1$ denote $\mathbf{y}_{l,p}=(y_{l},...,y_{p})$
and $\mathbf{S}%
_{l,p}^{(i)}=(s_{li},s_{l,i+1},s_{l+1,i},s_{l+1,i+1},...,s_{pi},s_{p,i+1}).$

Using (\ref{ExplX}) set%
\begin{equation*}
X_{i}^{(2)}\left( \mathbf{y}_{1,2};\mathbf{S}_{1,2}^{(i)}\right)
=X_{i}\left(
y_{1};s_{1i}X_{i}(y_{2};s_{2i},s_{2,i+1}),s_{1,i+1}s_{2,i+1}\right)
\end{equation*}%
and, by induction
\begin{equation*}
X_{i}^{(p)}\left( \mathbf{y}_{1,p};\mathbf{S}_{1,p}^{(i)}\right)
=X_{i}\left( y_{1};s_{1i}X_{i}^{(p-1)}\left( \mathbf{y}_{2,p};\mathbf{S}%
_{2,p}^{(i)}\right) ,\prod_{r=1}^{p}s_{r,i+1}\right) .
\end{equation*}%
Finally, recalling (\ref{ExplXstar}) put
\begin{equation*}
\bar{X}_{R_{i}}\left( \mathbf{y}_{1,p};\mathbf{S}_{1,p}^{(i)}\right)
=1-\left( 1-X_{i}^{(p)}\left( \mathbf{y}_{1,p};\mathbf{S}_{1,p}^{(i)}\right)
\right) ^{1/2^{i-1}}.
\end{equation*}%
It is not difficult to check that%
\begin{equation*}
\bar{X}_{R_{i}}\left( \mathbf{y}_{1,p};\mathbf{S}_{1,p}^{(i)}\right) =%
\mathbf{E}_{R_{i}}\left[
s_{1i}^{U_{ii}(y_{1})}s_{1,i+1}^{U_{i,i+1}(y_{1})}...s_{pi}^{U_{ii}(y_{p})}s_{p,i+1}^{U_{i,i+1}(y_{p})}%
\right] .
\end{equation*}

To complete the description of the limiting processes we are interesting in
introduce the function
\begin{equation*}
\psi (x;s)=\frac{1}{x+(1-x)/(1-s)},\quad s\in \lbrack 0,1],\,x\in \lbrack
0,1],
\end{equation*}%
and consider an $N-$dimensional process $\mathbf{U}_{N}(\cdot
)=(U_{N1}(\cdot ),...,U_{NN}(\cdot ))$ in which the first $N-1$ components
are equal to zero while $U_{NN}(\cdot )$ may be obtained by a time-change
from the following single-type continuous time Markov process $\sigma
(t),0\leq t<\infty $. The life-length distribution of particles in $\sigma
(\cdot )$ is exponential with parameter 1. Dying each particle produces
exactly two children. One may check (compare, for instance, with Example 3,
Section 8, Chapter 1 in \cite{Sev74}) that
\begin{equation*}
\mathbf{E}\left[ s^{\sigma (t)}|\sigma (0)=1\right] =1-\psi (1-e^{-t};s).
\end{equation*}%
Assuming that $\sigma (0)\overset{d}{=}R_{N}$ (recall (\ref{DefroN})) and
making the change of time $x=1-e^{-t},0\leq t<\infty ,$ we obtain an
inhomogeneous single-type branching process, denoted by $U_{NN}(\cdot )$
such that%
\begin{equation*}
\bar{G}_{R_{N}}(x;s)=\mathbf{E}_{R_{N}}\left[ s^{U_{NN}(x)}\right] =1-(\psi
(x;s))^{1/2^{N-1}}
\end{equation*}%
and%
\begin{equation*}
\mathbf{E}\left[ s^{U_{NN}(x+\Delta )}|U_{NN}(x)=1\right] =1-\psi \left(
\frac{\Delta }{1-x};s\right) ,\ 0<x+\Delta <1.
\end{equation*}%
Let, further, for $x_{j}\in \lbrack 0,1)$ and $\mathbf{S}%
_{j,p}=(s_{j},...,s_{p}),\,j=1,2,...,p$%
\begin{equation*}
G^{(1)}(x_{1};s_{1})=G(x_{1};s_{1})=1-\psi (x_{1};s_{1})
\end{equation*}%
and, by induction%
\begin{equation*}
G^{(p)}\left( \mathbf{x}_{1,p};\mathbf{S}_{1,p}\right) =G\left(
x_{1};s_{1N}G^{(p-1)}\left( \frac{\mathbf{x}_{2,p}}{1-x_{1}};\mathbf{S}%
_{2,p}\right) \right) .
\end{equation*}%
One may check that%
\begin{eqnarray*}
\bar{G}_{R_{N}}(\mathbf{x}_{1,p};\mathbf{S}_{1,p}) &=&\mathbf{E}_{R_{N}}%
\left[ s_{1}^{U_{NN}(x_{1})}s_{2}^{U_{NN}(x_{2})}...\,s_{p}^{U_{NN}(x_{p})}%
\right] \\
&=&1-(1-G^{(p)}(\mathbf{x}_{1,p};\mathbf{S}_{1,p}))^{1/2^{N-1}}.
\end{eqnarray*}

\section{Convergence of one-dimensional distributions\label{Sec4}}

As the first step in proving the main results of the paper we establish
convergence of one-dimensional distributions of $\{ \mathbf{Z}(m,n),0\leq
m\leq n\} $ given $\mathbf{Z}(n)\neq \mathbf{0}$. Let%
\begin{eqnarray*}
H_{m,n}^{\,(k,N)}(\mathbf{s})=\mathbf{E}\left[ \mathbf{s}^{\mathbf{Z}(m,n)}|%
\mathbf{Z}(0)=\mathbf{e}_{k}\right], J_{m,n}^{\,(k,N)}(\mathbf{s})=\mathbf{E}%
\left[ \mathbf{s}^{\mathbf{Z}(m,n)}|\mathbf{Z}(n)\neq \mathbf{0},\mathbf{Z}%
(0)=\mathbf{e}_{k}\right],
\end{eqnarray*}
\begin{eqnarray*}
\mathbf{H}_{m,n}(\mathbf{s}) =\left( H_{m,n}^{\,(1,N)}(\mathbf{s}%
),...,H_{m,n}^{\,(N,N)}(\mathbf{s})\right),\ \mathbf{J}_{m,n}(\mathbf{s}%
)=\left( J_{m,n}^{\,(1,N)}(\mathbf{s}),...,J_{m,n}^{\,(N,N)}(\mathbf{s}%
)\right).
\end{eqnarray*}%
For $\mathbf{x}=(x_{1},...,x_{N})$ and $\mathbf{y}=(y_{1},...,y_{N})$ put $%
\mathbf{x}\otimes \mathbf{y=}(x_{1}y_{1},x_{2}y_{2},...,x_{N}y_{N})$ and
denote
\begin{eqnarray}
s_{k}^{\prime }
&=&s_{k}Q_{n-m}^{(k,N)}+(1-Q_{n-m}^{(k,N)})=1-(1-s_{k})Q_{n-m}^{(k,N)},
\notag \\
&&  \notag \\
\mathbf{s}^{\prime } &=&(s_{1}^{\prime },...,s_{N}^{\prime })=\mathbf{1}-(%
\mathbf{1}-\mathbf{s})\otimes \mathbf{Q}_{n-m}.  \label{sprimeVect}
\end{eqnarray}%
It is not difficult to understand that
\begin{equation*}
H_{m,n}^{(k,N)}(\mathbf{s})=H_{m}^{(k,N)}(\mathbf{s}^{\prime
})=H_{m}^{(k,N)}(\mathbf{1}-(\mathbf{1}-\mathbf{s})\otimes \mathbf{Q}_{n-m})
\end{equation*}%
and that
\begin{equation}
J_{m,n}^{\,(k,N)}(\mathbf{s})=\mathbf{E}\left[ \mathbf{s}^{\mathbf{Z}(m,n)}|%
\mathbf{Z}(n)\neq \mathbf{0},\mathbf{Z}(0)=\mathbf{e}_{k}\right] =1-\frac{%
Q_{m}^{\,(k,N)}(\mathbf{s}^{\prime })}{Q_{n}^{\,(k,N)}}.  \label{BranProp2}
\end{equation}

\begin{theorem}
\label{T_familyMany}Let Hypothesis A be valid.

1) If $m\ll n^{\gamma _{1}}$ then
\begin{equation}
\lim_{n\rightarrow \infty }J_{m,n}^{\,(1,N)}(\mathbf{s})=\lim_{n\rightarrow
\infty }\mathbf{E}\left[ \mathbf{s}^{\mathbf{Z}(m,n)}|\mathbf{Z}(n)\neq
\mathbf{0},\mathbf{Z}(0)=\mathbf{e}_{1}\right] =s_{1}.  \label{LimSing}
\end{equation}%
2) If $n^{\gamma _{i}}\ll m\ll n^{\gamma _{i+1}}$ for some $i\in \left\{
1,2,...,N-1\right\} $ then
\begin{equation}
\lim_{n\rightarrow \infty }J_{m,n}^{\,(1,N)}(\mathbf{s}%
)=1-(1-s_{i+1})^{1/2^{i}}.  \label{LimSIngLate}
\end{equation}%
3) If $m=(y+l_{n})n^{\gamma _{i}},\,y\in \lbrack 0,\infty )$ for some $i\in
\left\{ 1,2,...,N-1\right\} $ then
\begin{equation}
\lim_{n\rightarrow \infty }J_{m,n}^{\,(1,N)}(\mathbf{s})=\bar{X}%
_{R_{i}}(y;s_{i},s_{i+1}).  \label{LimTWo}
\end{equation}%
4) If $m=(x+l_{n})n$, $x\in \lbrack 0,1)$ then
\begin{equation}
\lim_{n\rightarrow \infty }J_{m,n}^{{}\,(1,N)}(\mathbf{s})=\bar{G}%
_{R_{N}}(x;s_{N}).  \label{LimOLD}
\end{equation}
\end{theorem}

\textbf{Proof.} We start by observing that if $m\ll n$ then
\begin{eqnarray*}
1-s_{i}^{\prime } &=&(1-s_{i})Q_{n-m}^{(i,N)}\sim (1-s_{i})Q_{n}^{(i,N)} \\
&\sim &1-\exp \left\{ -(1-s_{i})Q_{n}^{(i,N)}\right\} \sim 1-\exp \left\{
-(1-s_{i})c_{iN}n^{-\gamma _{i}}\right\} .
\end{eqnarray*}%
This representation allows us to use the previous results with $s_{i}$ and $%
\lambda _{i}$ replaced by $s_{i}^{\prime }$ and $(1-s_{i})c_{iN}$,
respectively.

Recalling (\ref{SurvivSingle}) and applying Lemma \ref{L_OneOnly} we get%
\begin{equation*}
\lim_{n\rightarrow \infty }\frac{Q_{m}^{(1,N)}(\mathbf{s}^{\prime}) }{%
Q_{n}^{(1,N)}}=\lim_{n\rightarrow \infty }n^{2^{-(N-1)}}\frac{%
1-H_{m}^{(1,N)}( \mathbf{s}^{\prime }) }{C_{N}}=1-s_{1}.
\end{equation*}
Hence (\ref{LimSing}) follows.

Applying Lemma \ref{L_inbetween} with $n^{\gamma _{i}}\ll m\ll n^{\gamma
_{i+1}}$ and recalling Lemma \ref{L_DC} we conclude%
\begin{equation*}
\lim_{n\rightarrow \infty }\frac{Q_{m}^{(1,N)}(\mathbf{s}^{\prime })}{%
Q_{n}^{(1,N)}}=\frac{D_{i}}{C_{N}}%
((1-s_{i+1})c_{i+1,N})^{1/2^{i}}=(1-s_{i+1})^{1/2^{i}}
\end{equation*}%
leading to (\ref{LimSIngLate}).

\textbf{Proof of (\ref{LimTWo}).} If $y=0$ then the needed statement follows
from~(\ref{LimSing}) and~(\ref{LimSIngLate}). If $i\in \left\{
1,2,...,N-1\right\} $ is fixed and $m\sim yn^{\gamma _{i}},y>0,$ then for $%
j\geq i$
\begin{eqnarray*}
1-s_{j}^{\prime } &\sim &1-\exp \left\{ -(1-s_{j})c_{jN}n^{-\gamma
_{j}}\right\} \\
&\sim &1-\exp \left\{ -(1-s_{j})c_{jN}y^{2^{j-i}}m^{-2^{j-i}}\right\} .
\end{eqnarray*}%
Hence, by (\ref{SurvivSingle}) and Lemmas \ref{L_twoOnly} and \ref%
{L_MultiSharp} we get
\begin{equation*}
\lim_{n\rightarrow \infty }\frac{Q_{m}^{(1,N)}(\mathbf{s}^{\prime })}{%
Q_{n}^{(1,N)}}=\frac{D_{i-1}}{C_{N}}\left( \frac{\phi
_{i}(c_{iN}(1-s_{i})y,c_{i+1,N}(1-s_{i+1})y^{2})}{y}\right) ^{1/2^{i-1}}
\end{equation*}%
where we agree to write $D_{0}=1$. By~(\ref{DefSimpl}) and~(\ref{Const1})
\begin{eqnarray*}
&&\frac{\phi _{i}(c_{iN}(1-s_{i})y,c_{i+1,N}(1-s_{i+1})y^{2})}{y} \\
&&\quad =\sqrt{\frac{m_{i,i+1}c_{i+1,N}(1-s_{i+1})}{b_{i}}}\times \\
&&\qquad \times \frac{b_{i}c_{iN}(1-s_{i})+\sqrt{%
b_{i}m_{i,i+1}c_{i+1,N}(1-s_{i+1})}\tanh y\sqrt{%
b_{i}m_{i,i+1}c_{i+1,N}(1-s_{i+1})}}{b_{i}c_{iN}(1-s_{i})\tanh y\sqrt{%
b_{i}m_{i,i+1}c_{i+1,N}(1-s_{i+1})}+\sqrt{b_{i}m_{i,i+1}c_{i+1,N}(1-s_{i+1})}%
} \\
&&\quad =c_{iN}\sqrt{1-s_{i+1}}\times \frac{b_{i}c_{iN}(1-s_{i})+b_{i}c_{iN}%
\sqrt{1-s_{i+1}}\tanh (yb_{i}c_{iN}\sqrt{1-s_{i+1}})}{b_{i}c_{iN}(1-s_{i})%
\tanh (yb_{i}c_{iN}\sqrt{1-s_{i+1}})+b_{i}c_{iN}\sqrt{1-s_{i+1}}} \\
&&\quad =c_{iN}\sqrt{1-s_{i+1}}\times \frac{1-s_{i}+\sqrt{1-s_{i+1}}\tanh
(yb_{i}c_{iN}\sqrt{1-s_{i+1}})}{(1-s_{i})\tanh (yb_{i}c_{iN}\sqrt{1-s_{i+1}}%
)+\sqrt{1-s_{i+1}}}.
\end{eqnarray*}%
To complete the proof of (\ref{LimTWo}) it remains to recall (\ref{CD}).

\textbf{Proof of (\ref{LimOLD}).} If $x=0$ then (\ref{LimOLD}) follows from (%
\ref{LimSIngLate}). Consider now the case $m\sim xn,0<x<1$. Observe that for
$\mathbf{s}=(s_{1},s_{2},...,s_{N})\in \left[ 0,1\right] ^{N}$%
\begin{eqnarray}
H_{m}^{(1,N)}(\mathbf{1}^{(N-1)},s_{N})-H_{m}^{(1,N)}(\mathbf{s}) &=&\mathbf{%
E}\left[ \left( 1-s_{1}^{Z_{1}(m)}...\,s_{N-1}^{Z_{N-1}(m)}\right)
s_{N}^{Z_{N}(m)}\right]  \notag \\
\quad \leq \mathbf{E}\left[ 1-s_{1}^{Z_{1}(m)}...\,s_{N-1}^{Z_{N-1}(m)}%
\right] &\leq &\mathbf{P}(T_{N-1}>m)\leq cm^{-2^{-(N-2)}}.  \label{Neglig}
\end{eqnarray}%
Thus,%
\begin{equation*}
1-H_{m}^{(1,N)}(\mathbf{s})=1-H_{m}^{(1,N)}\left( \mathbf{1}%
^{(N-1)},s_{N}\right) +\varepsilon _{m,n}(\mathbf{s})Q_{m}^{(1,N)}
\end{equation*}%
where $\varepsilon _{m,n}(\mathbf{s})\rightarrow 0$ as $n\rightarrow \infty
, $ $m\sim xn$ $\ $uniformly in $\mathbf{s}\in \left[ 0,1\right] ^{N}.$
Therefore,
\begin{equation*}
1-H_{m}^{(1,N)}(\mathbf{s}^{\prime })=1-H_{m}^{(1,N)}\left( \hat{\mathbf{s}}%
,1-(1-s_{N})Q_{n-m}^{(N,N)}\right) +\varepsilon _{m,n}^{\prime }(\mathbf{s}%
)Q_{n}^{(1,N)}
\end{equation*}%
where $\varepsilon _{m,n}^{\prime }(\mathbf{s})\rightarrow 0$ as $%
n\rightarrow \infty ,$ $m\sim xn$ $\ $uniformly in $\hat{\mathbf{s}}=\left(
s_{1}^{\prime },...,s_{N-1}^{\prime }\right) \in \left[ 0,1\right] ^{N-1}$.

We now select an integer $r=r(m,n)\in \mathbb{N}^{\ast }=\left\{
1,2,...,\right\} $ in such a way that%
\begin{equation*}
H_{r-1}^{(N,N)}(0)\leq 1-(1-s_{N})Q_{n-m}^{(N,N)}\leq H_{r}^{(N,N)}(0)
\end{equation*}%
or%
\begin{equation*}
Q_{r}^{(N,N)}=1-H_{r}^{(N,N)}(0)\leq (1-s_{N})Q_{n-m}^{(N,N)}\leq
Q_{r-1}^{(N,N)}=1-H_{r-1}^{(N,N)}(0).
\end{equation*}%
This is possible, since by (\ref{SurvivSingle})
\begin{equation}
Q_{n-m}^{(N,N)}\sim \frac{1}{(n-m)b_{N}}\rightarrow 0,\ n-m\rightarrow
\infty .  \label{ASRep}
\end{equation}%
In particular,%
\begin{equation}
r\sim \frac{n-m}{1-s_{N}}.  \label{ASSR}
\end{equation}%
Under our choice of $r$, for any $\hat{\mathbf{s}}\in \left[ 0,1\right]
^{N-1}$
\begin{equation*}
H_{m}^{(1,N)}\left( \hat{\mathbf{s}},H_{r-1}^{(N,N)}(0)\right) \leq
H_{m}^{(1,N)}\left( \hat{\mathbf{s}},1-(1-s_{N})Q_{n-m}^{(N,N)}\right) \leq
H_{m}^{(1,N)}\left( \hat{\mathbf{s}},H_{r}^{(N,N)}(0)\right) .
\end{equation*}%
Letting $\hat{\mathbf{s}}=\left( H_{r}^{(1,N)}(\mathbf{0}%
),...,H_{r}^{(N-1,N)}(\mathbf{0})\right) $ we get by the branching property
of generating functions the estimate
\begin{equation*}
H_{m}^{(1,N)}\left( \hat{\mathbf{s}},1-(1-s_{N})Q_{n-m}^{(N,N)}\right) \leq
H_{m}^{(1,N)}(\mathbf{H}_{r}(\mathbf{0}))=H_{m+r}^{(1,N)}(\mathbf{0})
\end{equation*}%
implying
\begin{equation*}
1-H_{m}^{(1,N)}(\mathbf{s}^{\prime })\geq 1-H_{m+r}^{(1,N)}(\mathbf{0}%
)+\varepsilon _{m,n}^{\prime }Q_{n}^{(1,N)}=Q_{m+r}^{(1,N)}+\varepsilon
_{m,n}^{\prime }Q_{n}^{(1,N)},
\end{equation*}%
where $\varepsilon _{m,n}^{\prime }\rightarrow 0$ as $n\rightarrow \infty ,$
$m\sim xn,$ while $\hat{\mathbf{s}}=(H_{r-1}^{(1,N)}(\mathbf{0}%
),...,H_{r-1}^{(N-1,N)}(\mathbf{0}))$ gives the inequality
\begin{equation*}
H_{m}^{(1,N)}\left( \hat{\mathbf{s}},1-(1-s_{N})Q_{n-m}^{(N,N)}\right) \geq
H_{m}^{(1,N)}(\mathbf{H}_{r}(\mathbf{0}))=H_{m+r-1}^{(1,N)}(\mathbf{0})
\end{equation*}%
leading in the range under consideration to
\begin{equation*}
1-H_{m}^{(1,N)}(\mathbf{s}^{\prime })\leq Q_{m+r-1}^{(1,N)}+\varepsilon
_{m,n}^{\prime }Q_{n}^{(1,N)}.
\end{equation*}%
Hence%
\begin{equation*}
1-H_{m}^{(1,N)}(\mathbf{s}^{\prime })=Q_{m+r}^{(1,N)}+\varepsilon
_{m,n}^{\prime \prime }Q_{n}^{(1,N)}
\end{equation*}%
where $\varepsilon _{m,n}^{\prime \prime }\rightarrow 0$ as $n\rightarrow
\infty ,$ $m\sim xn$. We now conclude by (\ref{SurvivSingle}) that%
\begin{equation*}
1-H_{m}^{(1,N)}(\mathbf{s}^{\prime })\sim Q_{m+r}^{(1,N)}\sim
C_{N}(m+r)^{-2^{-(N-1)}}.
\end{equation*}%
Hence, on account of (\ref{ASSR}) and $m\sim xn,0<x<1,$ we get (recall (\ref%
{Const1}))%
\begin{eqnarray*}
\lim_{n\rightarrow \infty }\frac{1-H_{m}^{(1,N)}(\mathbf{s}^{\prime })}{%
Q_{n}^{(1,N)}} &=&\lim_{n\rightarrow \infty }\left( \frac{n}{%
nx+n(1-x)/(1-s_{N})}\right) ^{2^{-(N-1)}} \\
&=&\left( \frac{1}{x+(1-x)/(1-s_{N})}\right) ^{2^{-(N-1)}}
\end{eqnarray*}%
completing the proof of (\ref{LimOLD}).

Theorem \ref{T_familyMany} is proved.\

\textbf{Proof} \textbf{of Theorem \ref{T_Yaglom}}. Since our process is
decomposable and strongly critical, it is sufficient to check (\ref{Yag})
for $i=1$ only. For $\hat{s}_{N}=\exp (-\lambda /(nb_{N}))$ we have%
\begin{equation*}
\mathbf{E}\left[ \exp \left\{ -\lambda \frac{Z_{N}(n)}{b_{N}n}\right\} \Big|%
\mathbf{Z}(n)\neq \mathbf{0}\right] =1-\frac{1-H_{n}^{(1,N)}(\mathbf{1}%
^{(N-1)},\hat{s}_{N})}{Q_{n}^{(1,N)}}.
\end{equation*}%
We now select an integer $r=r(\lambda ,n)\in \mathbb{N}^{\ast }=\left\{
1,2,...,\right\} $ in such a way that%
\begin{equation*}
H_{r-1}^{(N,N)}(0)\leq \hat{s}_{N}\leq H_{r}^{(N,N)}(0).
\end{equation*}%
It follows from (\ref{ASRep})\ that $r\sim n\lambda ^{-1}$. Letting $%
s_{i}=H_{r}^{(i,N)}(\mathbf{0}),i=1,2,...,N,$ and setting $\mathbf{s}%
=(s_{1},...,s_{N})$ we get by (\ref{Neglig}) after evident estimates that%
\begin{equation*}
\left\vert H_{n}^{(1,N)}(\mathbf{1}^{(N-1)},\hat{s}_{N})-H_{n}^{(1,N)}(%
\mathbf{s})\right\vert \leq cn^{1/2^{N-2}}.
\end{equation*}%
Hence, using (\ref{SurvivSingle}) with $i=1$ we obtain
\begin{eqnarray*}
\frac{1-H_{n}^{(1,N)}(\mathbf{1}^{(N-1)},\hat{s}_{N})}{Q_{n}^{(1,N)}} &\sim &%
\frac{1-H_{n}^{(1,N)}(\mathbf{s})}{Q_{n}^{(1,N)}}=\frac{Q_{r+n}^{(1,N)}}{%
Q_{n}^{(1,N)}} \\
&\sim &\left( \frac{n}{r+n}\right) ^{1/2^{N-1}}\sim \left( \frac{\lambda }{%
1+\lambda }\right) ^{1/2^{N-1}}
\end{eqnarray*}%
as desired.\

\section{Convergence of finite-dimensional distributions\label{Sec5}}

In this section we study the limiting behavior of the finite-dimensional
distributions of the reduced process $\left\{ \mathbf{Z}( m,n) ,0\leq m\leq
n\right\} $. Our first theorem deals with the case $m\ll n$.

\begin{theorem}
\label{T_findim} Let Hypothesis A be valid and $\mathbf{S}%
_{l}=(s_{l1},...,s_{lN}),l=1,2,...,p$.

1) If, for a fixed $i\in \left\{ 0,1,...,N-1\right\}$
\begin{equation*}
n^{\gamma _{i}}\ll m_{l}\ll n^{\gamma _{i+1}},\quad l=1,...,p
\end{equation*}%
then%
\begin{equation}
\lim_{n\rightarrow \infty }\mathbf{E}\left[ \prod_{l=1}^{p}\mathbf{S}_{l}^{%
\mathbf{Z}(m_{l},n)}\,\Big|\,\mathbf{Z}(n)\neq 0\right] =1-\left(
1-\prod_{l=1}^{p}s_{l,i+1}\right) ^{1/2^{i}}.  \label{Lim_diskr}
\end{equation}

2) Let $0=Y_{1}<Y_{2}<...<Y_{p}<\infty $ be a tuple of nonnegative numbers
with $y_{1}=0,$ $y_{l}=Y_{l}-Y_{l-1},\,l=2,...,p.$ If, for a fixed $i\in
\left\{ 1,2,...,N-1\right\} $
\begin{equation*}
m_{1}\sim l_{n}n^{\gamma _{i}},\quad m_{l}\sim Y_{l}n^{\gamma
_{i}},\,l=2,...,p
\end{equation*}%
then
\begin{equation}
\lim_{n\rightarrow \infty }\mathbf{E}\left[ \prod_{l=1}^{p}\mathbf{S}_{l}^{%
\mathbf{Z}(m_{l},n)}\,\Big|\,\mathbf{Z}(n)\neq 0\right] =\bar{X}_{R_{i}}\Big(%
\mathbf{y}_{1,p};\mathbf{S}_{1,p}^{(i)}\Big).  \label{Lim_Multi}
\end{equation}
\end{theorem}

The second theorem is devoted to the finite-dimensional distributions of the
reduced process when $m$ is of order $n$.

\begin{theorem}
\label{T_endpoint}Let\textbf{\ }Hypothesis A be valid and $%
0=X_{1}<X_{2}<...<X_{p}<1$ be a tuple of nonnegative numbers with $x_{1}=0,$
$x_{l}=X_{l}-X_{l-1},\,l=2,...,p$. If
\begin{equation*}
m_{1}\sim l_{n}n,\quad m_{l}\sim X_{l}n,\,l=2,...,p
\end{equation*}%
then
\begin{equation*}
\lim_{n\rightarrow \infty }\mathbf{E}\left[ \prod_{l=1}^{p}\mathbf{S}_{l}^{%
\mathbf{Z}(m_{l},n)}\,\Big|\,\mathbf{Z}(n)\neq 0\right] =\bar{G}_{R_{N}}\Big(%
\mathbf{x}_{1,p};\mathbf{S}_{1,p;N}\Big),
\end{equation*}%
where $\mathbf{S}_{1,p;N}=(s_{1N},s_{2N},...,s_{pN})$.
\end{theorem}

To prove Theorems \ref{T_findim} and \ref{T_endpoint} we need additional
notation.

For $0\leq m_{0}<m_{1}<...<m_{p}\leq n$ set $\mathbf{m}%
=(m_{0},m_{1},...,m_{p}),$ put $\Delta _{i}=m_{i}-m_{i-1}$, and denote%
\begin{eqnarray*}
\hat{J}_{m_{0},m_{1},...,m_{p},n}^{\,(i,N)}(\mathbf{S}_{1},...,\mathbf{S}%
_{p}) &=&\hat{J}_{\mathbf{m},n}^{\,(i,N)}(\mathbf{S}_{1},...,\mathbf{S}_{p})
\\
&=&\mathbf{E}\left[ \prod_{l=1}^{p}\mathbf{S}_{l}^{\mathbf{Z}(m_{l},n)}\,%
\Big|\,\mathbf{Z}(m_{0},n)=\mathbf{e}_{i}\right]
\end{eqnarray*}%
and%
\begin{equation*}
\mathbf{\hat{J}}_{\mathbf{m},n}(\mathbf{S}_{1},...,\mathbf{S}_{p})=\left(
\hat{J}_{\mathbf{m},n}^{\,(1,N)}(\mathbf{S}_{1},...,\mathbf{S}_{p}),...,\hat{%
J}_{\mathbf{m},n}^{\,(N,N)}(\mathbf{S}_{1},...,\mathbf{S}_{p})\right) .
\end{equation*}

The next statement is a simple observation following from Corollary 2 in~%
\cite{VD06}.

\begin{lemma}
\label{L_convol}For any $0\leq m_{0}<m_{1}<...<m_{p}\leq n$ we have%
\begin{eqnarray*}
&&\hat{J}_{\mathbf{m},n}^{\,\,(1,N)}(\mathbf{S}_{1},...,\mathbf{S}_{p})=\hat{%
J}_{m_{0},m_{1},n}^{\,(1,N)}\left( \mathbf{S}_{1}\otimes \mathbf{\hat{J}}%
_{m_{1},m_{2},...,m_{p},n}(\mathbf{S}_{2},...,\mathbf{S}_{p})\right) \\
&&\quad=J_{\Delta _{1},n-m_{0}}^{\,(1,N)}\left( \mathbf{S}_{1}\otimes
\mathbf{J}_{\Delta _{2},n-m_{1}}\left( \mathbf{S}_{2}\otimes ...(\mathbf{S}%
_{p-1}\otimes \mathbf{J}_{\Delta _{p},n-m_{p-1}}(\mathbf{S}_{p})\right)
....)\right) .
\end{eqnarray*}%
In particular, if $\mathbf{m}=(0,m_{1},m_{2})$ then%
\begin{equation*}
\hat{J}_{\mathbf{m},n}^{\,(1,N)}(\mathbf{S}_{1},\mathbf{S}%
_{2})=J_{m,n}^{\,(1,N)}(\mathbf{S}_{1}\otimes \mathbf{J}_{\Delta
_{2},n-m_{1}}(\mathbf{S}_{2}))
\end{equation*}%
and if $\mathbf{m}=(m_{0},m_{1})$ then for $\mathbf{s}=(s_{1},...,s_{N})$%
\begin{equation}
\hat{J}_{m_{0},m_{1},n}^{\,(k,N)}(\mathbf{s})=J_{\Delta
_{1},n-m_{0}}^{\,(k,N)}(\mathbf{s})=1-\frac{1-H_{\Delta _{1}}^{(k,N)}\left(
\mathbf{1}-(\mathbf{1}-\mathbf{s})\otimes \mathbf{Q}_{n-m_{1}}\right) }{%
Q_{n-m_{0}}^{(k,N)}}.  \label{Derivat}
\end{equation}
\end{lemma}

Using (\ref{Derivat}) we prove the following statement.

\begin{lemma}
\label{L_derivative}If $m_{0}=(Y_{0}+l_{n})n^{\gamma
_{i}}<m_{1}=(Y_{1}+l_{n})n^{\gamma _{i}}$ then for any $j\geq i$ there
exists a constant $\chi \in (0,\infty )$ such that for all $n\geq n_{0}$ \
\begin{equation*}
\mathbf{P}(\mathbf{Z}(m_{1},n)=\mathbf{e}_{j}|\mathbf{Z}(m_{0},n)=\mathbf{e}%
_{j})\geq 1-\chi (Y_{1}-Y_{0}).
\end{equation*}
\end{lemma}

\textbf{Proof.} By the decomposability assumption and the condition $%
m_{jj}=1 $ implying
\begin{equation*}
m_{jj}(\Delta _{1})=\frac{\partial H_{\Delta _{1}}^{(j,N)}(\mathbf{s})}{%
\partial s_{j}}\left\vert _{\mathbf{s}=\mathbf{1}}\right. =1
\end{equation*}%
we get
\begin{equation*}
1-\frac{\partial H_{\Delta _{1}}^{(j,N)}(\mathbf{s})}{\partial s_{j}}%
\left\vert _{\mathbf{s}=\mathbf{H}_{n-m_{1}}(\mathbf{0})}\right. \leq
\sum_{k=j}^{N}\mathbf{E}Z_{j}(\Delta _{1})(Z_{k}(\Delta _{1})-\delta
_{kj})Q_{n-m_{1}}^{(k,N)}.
\end{equation*}%
Recalling (\ref{Momvariance}) and (\ref{SurvivSingle}) and setting $%
h=Y_{1}-Y_{0}$ we obtain
\begin{eqnarray*}
\mathbf{E}Z_{j}(\Delta _{1})Z_{k}(\Delta _{1})Q_{n-m_{1}}^{(k,N)} &\leq
&c_{0}(n-m_{1})^{-1/2^{N-k}}(\Delta _{1})^{k-j+1} \\
&\leq &c_{0}(n-m_{1})^{-1/2^{N-k}}(hn^{1/2^{N-i}})^{k-j+1} \\
&\leq &\chi hn^{-1/2^{N-k}}(n^{1/2^{N-i}})^{k-j+1}
\end{eqnarray*}%
for some constants $0<c_{0}\leq \chi <\infty $. On account of $k\geq j\geq i$
we have%
\begin{equation*}
\frac{k-j+1}{2^{N-i}}-\frac{1}{2^{N-k}}=\frac{1}{2^{N-i}}(k-j+1-2^{k-i})\leq
0.
\end{equation*}%
Thus,%
\begin{equation*}
1-\frac{\partial H_{\Delta _{1}}^{(j,N)}(\mathbf{s})}{\partial s_{j}}%
\left\vert _{\mathbf{s}=\mathbf{H}_{n-m_{1}}(\mathbf{0})}\right. \leq \chi h.
\end{equation*}%
Hence, using the previous lemma and monotonicity of $Q_{r}^{(j,N)}$ in $r$
we get%
\begin{eqnarray}
\mathbf{P}(\mathbf{Z}(m_{1},n)=\mathbf{e}_{j}|\mathbf{Z}(m_{0},n)=\mathbf{e}%
_{j}) &=&\frac{Q_{n-m_{1}}^{(j,N)}}{Q_{n-m_{0}}^{(j,N)}}\frac{\partial
H_{\Delta _{1}}^{(j,N)}(\mathbf{s})}{\partial s_{j}}\left\vert _{\mathbf{s}=%
\mathbf{H}_{n-m_{1}}(\mathbf{0})}\right.  \label{ExplCoeff} \\
&\geq &\frac{\partial H_{\Delta _{1}}^{(j,N)}(\mathbf{s})}{\partial s_{j}}%
\left\vert _{\mathbf{s}=\mathbf{H}_{n-m_{1}}(\mathbf{0})}\right. \geq 1-\chi
h.  \notag
\end{eqnarray}%
Lemma \ref{L_derivative} is proved. \

\textbf{Proof} \textbf{of Theorem \ref{T_findim}}. Using (\ref{BranProp2})
and Theorem \ref{T_familyMany} we see that

1) if $m\ll n^{\gamma _{k}}$ then%
\begin{equation}
\lim_{n\rightarrow \infty }J_{m,n}^{\,(k,N)}( \mathbf{s}) =s_{k};
\label{LeftEnd}
\end{equation}

2) if $m=(y+l_{n})n^{\gamma _{k}}=(y+l_{n})n^{1/2^{( N-k) }},\,y\in \lbrack
0,\infty )$ then
\begin{equation*}
\lim_{n\rightarrow \infty }J_{m,n}^{\,(k,N)}( \mathbf{s}) =X_{k}(
y;s_{k},s_{k+1}) ;
\end{equation*}

3) if $m=(x+l_{n})n,\,x\in \left[ 0,1\right] $ then%
\begin{equation}
\lim_{n\rightarrow \infty }J_{m,n}^{\,(N,N)}(\mathbf{s})=G\left(
x;s_{N}\right) .  \label{End1}
\end{equation}

\textbf{Proof of (\ref{Lim_diskr})}. Consider first the case $p=2$ and take $%
\mathbf{m}=(0,m_{1},m_{2})$. By Lemma \ref{L_convol}%
\begin{equation}
\hat{J}_{\mathbf{m},n}^{\,(1,N)}(\mathbf{S}_{1},\mathbf{S}%
_{2})=J_{m_{1},n}^{\,(1,N)}(\mathbf{S}_{1}\otimes \mathbf{J}_{\Delta
_{2},n-m_{1}}(\mathbf{S}_{2})).  \label{JointDist}
\end{equation}

It follows from (\ref{LimSIngLate}) that, given $n^{\gamma _{i}}\ll m_{1}\ll
n^{\gamma _{i+1}}$
\begin{equation*}
J_{m_{1},n}^{\,(1,N)}( \mathbf{S}_{1}) \rightarrow 1-( 1-s_{1,i+1})
^{1/2^{i}}
\end{equation*}%
as $n\rightarrow \infty.$ Further, in view of $\Delta _{2}=m_{2}-m_{1}\ll
n^{\gamma _{i+1}}$ and~(\ref{LeftEnd}) $J_{\Delta _{2},n-m_{1}}^{\,(i+1,N)}(%
\mathbf{S}_{2})\rightarrow s_{2,i+1}$ as $n\rightarrow \infty.$ Hence, using
the continuity of the functions under consideration and (\ref{JointDist}) we
get%
\begin{equation*}
\lim_{n\rightarrow \infty }\hat{J}_{\mathbf{m},n}^{\,(1,N)}(\mathbf{S}_{1},%
\mathbf{S}_{2})=1-( 1-s_{1,i+1}s_{2,i+1}) ^{1/2^{i}}.
\end{equation*}

The validity of (\ref{Lim_diskr}) for any $p>3$ may be checked by induction
using Lemma~\ref{L_convol}.

\textbf{Proof of (\ref{Lim_Multi})}. Consider again the case $p=2$ only. It
follows from (\ref{LimTWo}) that, given $m_{l}\sim Y_{l}n^{\gamma
_{i}},\,l=1,2,$ with $Y_{1}=y_{1}$
\begin{equation*}
J_{m_{1},n}^{\,(1,N)}(\mathbf{s})\rightarrow \bar{X}%
_{R_{i}}(y_{1};s_{i},s_{i+1})
\end{equation*}%
as $n\rightarrow \infty $ and%
\begin{equation*}
\lim_{n\rightarrow \infty }J_{\Delta _{2},n-m_{1}}^{\,(i,N)}(\mathbf{S}%
_{2})=X_{i}(y_{2};s_{2i},s_{2,i+1}),\quad \lim_{n\rightarrow \infty
}J_{\Delta _{2},n-m_{1}}^{\,(i+1,N)}(\mathbf{S}_{2})=s_{2,i+1}.
\end{equation*}%
Hence, using the continuity of the functions involved and (\ref{JointDist})
we get%
\begin{equation*}
\lim_{n\rightarrow \infty }\hat{J}_{\mathbf{m},n}^{\,(1,N)}(\mathbf{S}_{1},%
\mathbf{S}_{2})=\bar{X}_{R_{i}}\left(
y_{1};s_{1i}X_{i}(y_{2};s_{2i},s_{2,i+1}),s_{1,i+1}s_{2,i+1}\right)
\end{equation*}%
proving (\ref{Lim_Multi}) for $p=2$.

To justify (\ref{Lim_Multi}) for $p>3$ it is necessary to use Lemma \ref%
{L_convol} and induction arguments. We omit the respective details. \

\textbf{Proof} \textbf{of Theorem \ref{T_endpoint}}. We consider the case $%
p=2$ only and to this aim take $\mathbf{m}%
=(0,(x_{1}+l_{n})n,(x_{1}+x_{2}+l_{n})n)$. By (\ref{JointDist}), (\ref%
{LimOLD}) and (\ref{End1})
\begin{eqnarray*}
\lim_{n\rightarrow \infty }\hat{J}_{\mathbf{m},n}^{\,(1,N)}(\mathbf{S}_{1},%
\mathbf{S}_{2}) &=&\lim_{n\rightarrow \infty }J_{(x_{1}+l_{n})n,n}^{\,(1,N)}(%
\mathbf{S}_{1}\otimes \mathbf{J}_{x_{2}n,n(1-x_{1}-l_{n})}(\mathbf{S}_{2}))
\\
&=&\bar{G}_{R_{N}}\left( x_{1};s_{1N}G\left( \frac{x_{2}}{1-x_{1}}%
;s_{2N}\right) \right) =\bar{G}_{R_{N}}(\mathbf{x}_{1,2};\mathbf{S}_{1,2;N}).
\end{eqnarray*}%
The desired statement for $p>2$ follows by induction. \

\textbf{Proof} \textbf{of point 1) of Theorem \ref{T_Skhod1}}. Let $%
0=t_{0}<t_{1}<...<t_{p}<1$. If $\gamma _{i-1}\leq t_{1}<t_{p}<\gamma _{i}$
for some $i\in \left\{ 1,2,...,N\right\} $ then the needed convergence of
finite-dimensional distributions follows from (\ref{Lim_diskr}). We now
consider another extreme case, namely, take a tuple $%
0=t_{0}<t_{1}<...<t_{N}<1$ such that $\gamma _{i-1}\leq t_{i}<\gamma _{i}$
for all $i=1,2,...,N$. Then for $m_{i}\sim n^{t_{i}}g_{n}(t_{i})$ we have
\begin{equation*}
n^{\gamma _{i-1}}\ll m_{i}\ll n^{\gamma _{i}},\quad \Delta
_{i}=m_{i}-m_{i-1}\sim m_{i},\quad n-m_{i}\sim n.
\end{equation*}%
These relations, (\ref{LimSIngLate}), (\ref{LeftEnd}), and the continuity of
the respective probability generating functions imply (recall (\ref{DefOmega}%
))%
\begin{eqnarray*}
&&\lim_{n\rightarrow \infty }\hat{J}_{\mathbf{m},n}^{\,(1,N)}(\mathbf{S}%
_{1},...,\mathbf{S}_{N}) \\
&&\quad =\lim_{n\rightarrow \infty }J_{m_{1},n}^{\,(1,N)}(\mathbf{S}%
_{1}\otimes \mathbf{J}_{m_{2},n}(\mathbf{S}_{2}\otimes ...(\mathbf{S}%
_{N-1}\otimes \mathbf{J}_{m_{N},n}(\mathbf{S}_{N}))...)) \\
&&\quad =s_{11}\left( 1-\sqrt{1-\lim_{n\rightarrow \infty
}J_{m_{2},n}^{\,(2,N)}(\mathbf{S}_{2}\otimes ...(\mathbf{S}_{N-1}\otimes
\mathbf{J}_{m_{N},n}(\mathbf{S}_{N}))...)}\,\right) \\
&&\quad =s_{11}\left( 1-\sqrt{1-s_{22}\left( 1-\sqrt{1-\Omega
_{N-2}(s_{33},...,s_{NN})}\,\right) }\,\,\right) \\
&&\quad =...=\Omega _{N}(s_{11},s_{22},...,s_{NN})
\end{eqnarray*}%
as required.

The case when several values among $t_{j}$ are contained in a subinterval $%
[\gamma _{i-1},\gamma _{i})$ may be considered by combining the previous
arguments. We omit the respective details. \

\section{Tightness\label{Sec6}}

Denote by $\mathbf{z}^{(i,i+1)},1\leq i\leq N-1,$ the $(N-2)$-dimensional
vector obtained from $\mathbf{z}=(z_{1},...,z_{N})\in \mathbb{Z}_{+}^{N}$ by
deleting the coordinates $i$ and $i+1$ and by $\mathbf{z}^{(i)},1\leq i\leq
N-1,$ the $(N-1)$-dimensional vector obtained from $\mathbf{z}$ by deleting
the $i$-th coordinate. Let $\left\Vert \mathbf{x}\right\Vert $ be the sum of
absolute values of all coordinates of the vector $\mathbf{x}$.

Set $\mathcal{C}_{i}=\left\{ \mathbf{z}\in \mathbb{Z}_{+}^{N}:\left\Vert
\mathbf{z}^{(i)}\right\Vert >0\right\} ,\,\mathcal{B}_{i}=\mathbb{Z}%
_{+}^{N}\backslash \mathcal{C}_{i}$ and
\begin{equation*}
\mathcal{C}_{i,i+1}=\left\{ \mathbf{z}\in \mathbb{Z}_{+}^{N}:\left\Vert
\mathbf{z}^{(i,i+1)}\right\Vert >0\right\} .
\end{equation*}

Put \b{Z}$_{i}(m)=Z_{1}(m)+...+Z_{i}(m)$ and denote%
\begin{equation*}
\text{\b{Z}}_{i}(m,n)=\sum_{k=1}^{i}Z_{k}(m,n),\quad \bar{Z}%
_{i}(m,n)=\sum_{k=i}^{N}Z_{k}(m,n).
\end{equation*}

In what follows it will be convenient to write $\mathbf{P}_{n}(\mathcal{B})$
for $\mathbf{P}(\mathcal{B}|\mathbf{Z}(n)\neq \mathbf{0},\mathbf{Z}(0)=%
\mathbf{e}_{1})$ for any admissible event~$\mathcal{B}$.

We start checking the desired tightness of the prelimiting processes in
Theorems \ref{T_SkorohConst} and~\ref{T_Skhod1} by proving two important
lemmas.

Let $A_{i}(n) =\left\{ m:n^{\gamma _{i}}g_{n}(\gamma _{i})\leq m<n^{\gamma
_{i+1}-\varepsilon }g_{n}(\gamma _{i+1}-\varepsilon
)\right\},\,\varepsilon>0.$

\begin{lemma}
\label{L_Negli} For any $i=0,1,2,...,N-1$ and $\varepsilon\in (0,\gamma_1)$
\begin{equation*}
\lim_{n\rightarrow \infty }\mathbf{P}_{n}( \exists m\in A_{i}( n) :\mathbf{Z}%
( m,n) \in \mathcal{C}_{i+1}) =0.
\end{equation*}
\end{lemma}

\textbf{Proof.} If $m\in A_{i}(n)$ then $\bar{Z}_{i+2}(m,n)\leq \bar{Z}%
_{i+2}(n^{\gamma _{i+1}-\varepsilon }g_{n}(\gamma _{i+1}-\varepsilon ),n)$
and
\begin{equation*}
\left\{ \text{\b{Z}}_{i}(m,n)>0\right\} \Rightarrow \left\{ \text{\b{Z}}%
_{i}(m)>0\right\} \Rightarrow \left\{ \text{\b{Z}}_{i}(n^{\gamma
_{i}}g_{n}(\gamma _{i}))>0\right\} .
\end{equation*}%
Thus,
\begin{eqnarray*}
\mathbf{P}_{n}(\exists m\in A_{i}(n):\mathbf{Z}(m,n)\in \mathcal{C}_{i+1})
&\leq &\mathbf{P}_{n}(\bar{Z}_{i+2}(n^{\gamma _{i+1}-\varepsilon
}g_{n}(\gamma _{i+1}-\varepsilon ),n)>0) \\
&&+\mathbf{P}_{n}\left( \text{\b{Z}}_{i}(n^{\gamma _{i}}g_{n}(\gamma
_{i}))>0\right) .
\end{eqnarray*}%
Letting $n$ tend to infinity we see that the first summand at the right-hand
side of the inequality vanishes by (\ref{LimSIngLate}), while the second one
is zero for $i=0$ and tends to zero for $1\leq i\leq N-1$ in view of
\begin{eqnarray*}
\mathbf{P}_{n}\left( \text{\b{Z}}_{i-1}(n^{\gamma _{i}}g_{n}(\gamma
_{i}))>0\right) &=&\frac{\mathbf{P}(T_{i}>n^{\gamma _{i}}g_{n}(\gamma _{i}))%
}{\mathbf{P}(T_{N}>n)} \\
&\sim &\frac{c_{1i}}{c_{1N}}\frac{n^{1/2^{N-1}}}{(n^{1/2^{N-i}}g_{n}(\gamma
_{i}))^{1/2^{i-1}}}=\frac{c_{1i}}{c_{1N}}\frac{1}{(g_{n}(\gamma
_{i}))^{1/2^{i-1}}}.
\end{eqnarray*}%
The lemma is proved.\

\begin{lemma}
\label{L_Skor1}If $N\geq 3$ then for any $i=1,2,...,N-1$%
\begin{equation*}
\lim_{n\rightarrow \infty }\mathbf{P}_{n}( \exists m\in \lbrack n^{3\gamma
_{i-1}},n^{3\gamma _{i}}]:\mathbf{Z}( m,n) \in \mathcal{C}_{i,i+1}) =0.
\end{equation*}
\end{lemma}

\textbf{Proof.} By the same arguments as in Lemma \ref{L_Negli}, we conclude
\begin{eqnarray*}
&&\mathbf{P}_{n}(\exists m\in \lbrack n^{3\gamma _{i-1}},n^{3\gamma _{i}}]:%
\mathbf{Z}(m,n)\in \mathcal{C}_{i,i+1}) \\
&&\qquad \leq \mathbf{P}_{n}(\bar{Z}_{i+2}(n^{3\gamma _{i}},n)>0)+\mathbf{P}%
_{n}(\text{\b{Z}}_{i-1}(n^{3\gamma _{i-1}})>0).
\end{eqnarray*}%
According to point 3) of Theorem \ref{T_familyMany} the first summand tends
to zero as $n\rightarrow \infty $ while the second is, by definition zero
for $i=1$ and is evaluated as%
\begin{equation*}
\frac{\mathbf{P}(T_{i-1}>n^{3\gamma _{i-1}})}{\mathbf{P}(T_{N}>n)}\sim \frac{%
c_{1,i-1}}{c_{1N}}\frac{n^{1/2^{N-1}}}{(n^{3/2^{N-i+1}})^{1/2^{i-2}}}\sim
\frac{c_{1,i-1}}{c_{1N}}\frac{1}{n^{1/2^{N-2}}}
\end{equation*}%
for $i\geq 2$. \ This completes the proof of the lemma.

\subsection{Macroscopic view\label{Sec51}}

In this section we prove Theorem \ref{T_SkorohConst}\textbf{\ }which
describes the macroscopic structure of the family tree. Convergence of the
finite-dimensional distributions of $\left\{ \mathbf{Z}(n^{t}g_{n}(t),n),0%
\leq t<1\right\} $ to the respective finite-dimensional distributions of $%
\left\{ \mathbf{R}(t),0\leq t<1\right\} $ has been established in (\ref%
{Lim_diskr}). Thus, we concentrate on proving the tightness.

Since $\mathbf{Z}(n^{t}g_{n}(t),n)$ has integer-valued components we need to
check\ for each interval $A_{i}=\left[ \gamma _{i},\gamma _{i+1}-\varepsilon %
\right] ,i=0,1,...,N-1,$ that (see \cite{Bil68}, Theorem 15.3)

1) for any positive $\eta $ there exists $L$ such that
\begin{equation}
\mathbf{P}_{n}\left( \sup_{t\in A_{i}}\left\Vert \mathbf{Z}%
(n^{t}g_{n}(t),n)\right\Vert >L\right) \leq \eta ,~n\geq 1;  \label{Bi0}
\end{equation}

2) for any positive $\eta $ there exist $\delta >0$ and $n_{0}$ such that,
for all $n\geq n_{0}$%
\begin{equation}
\mathbf{P}_{n}\left( \max \left( \min_{k=1,2}\left\Vert \mathbf{Z}%
(n^{t}g_{n}(t),n)-\mathbf{Z}(n^{t_{k}}g_{n}(t_{k}),n)\right\Vert \right)
\neq 0\right) \leq \eta ,  \label{Bi01}
\end{equation}%
where the $\max $ is taken over all $\gamma _{i}\leq t_{1}\leq t\leq
t_{2}\leq \gamma _{i+1}-\varepsilon $ such that $t_{2}-t_{1}\leq \delta ;$

\begin{equation}
\mathbf{P}_{n}(\exists t,s\in \left[ \gamma _{i},\gamma _{i}+\delta \right] :%
\mathbf{Z}(n^{t}g_{n}(t),n)\neq \mathbf{Z}(n^{s}g_{n}(s),n)\,)\leq \eta ,
\label{Bi02}
\end{equation}%
and%
\begin{equation}
\mathbf{P}_{n}(\exists t,s\in \lbrack \gamma _{i+1}-\delta -\varepsilon
,\gamma _{i+1}-\varepsilon ]:\mathbf{Z}(n^{t}g_{n}(t),n)\neq \mathbf{Z}%
(n^{s}g_{n}(s),n)\,\,)\leq \eta .  \label{Bi03}
\end{equation}

The fact that the random variable $\left\Vert \mathbf{Z}( n^{t}g_{n}(t),n)
\right\Vert $ is monotone in $t$ for fixed $n$ essentially simplifies the
proof.

Indeed, in this case
\begin{equation*}
\mathbf{P}_{n}\left( \sup_{t\in A_{i}}\left\Vert \mathbf{Z}%
(n^{t}g_{n}(t),n)\right\Vert >L\right) \leq \mathbf{P}_{n}(\left\Vert
\mathbf{Z}(n^{1-\varepsilon }g_{n}(1-\varepsilon ),n)\right\Vert >L\,)
\end{equation*}%
and (\ref{Bi0}) follows from the one-dimensional convergence established in (%
\ref{LimSIngLate}) for $i=N-1$.

To prove (\ref{Bi01})-(\ref{Bi03}) we introduce the events%
\begin{eqnarray*}
\mathcal{D}_{i} &=&\left\{ \forall t\in A_{i}:\mathbf{Z}(n^{t}g_{n}(t),n)\in
\mathcal{B}_{i+1}\right\} , \\
\mathcal{F}_{i}(a,b) &=&\left\{ \exists t,s\in \lbrack
a,b]:Z_{i+1}(n^{t}g_{n}(t),n)\neq Z_{i+1}(n^{s}g_{n}(s),n)\right\} ,
\end{eqnarray*}%
take a sufficiently small $\delta >0$ and observe that if $\left[ a,b\right]
\subset \left[ \gamma _{i},\gamma _{i+1}-\varepsilon \right] $ then
\begin{eqnarray*}
&&\mathbf{P}_{n}(\exists t,s\in \lbrack a,b]:\mathbf{Z}(n^{t}g_{n}(t),n)\neq
\mathbf{Z}(n^{s}g_{n}(s),n)\,) \\
&&\qquad \leq \mathbf{P}_{n}(\exists t\in A_{i}:\mathbf{Z}%
(n^{t}g_{n}(t),n)\in \mathcal{C}_{i+1})+\mathbf{P}_{n}(\mathcal{D}_{i}\cap
\mathcal{F}_{i}(\gamma _{i},\gamma _{i+1}-\varepsilon )).
\end{eqnarray*}%
By Lemma \ref{L_Negli} the first term at the right-hand side tends to zero
as $n\rightarrow \infty $.

Further, for $i\geq 1$
\begin{equation*}
\mathbf{P}_{n}( \mathcal{D}_{i}\cap \mathcal{F}_{i}( \gamma _{i},\gamma
_{i+1}-\varepsilon ) ) \leq \mathbf{P}_{n}( Z_{i+1}( n^{\gamma _{i}}g_{n},n)
\neq Z_{i+1}( n^{\gamma _{i+1}-\varepsilon }g_{n},n) ) \rightarrow 0
\end{equation*}%
by (\ref{Lim_diskr}). This justifies (\ref{Bi02})-(\ref{Bi03}).

To check the validity of (\ref{Bi01}) it remains to note that

\begin{eqnarray*}
&&\mathbf{P}_{n}\left( \max \left( \min_{k=1,2}\left\Vert \mathbf{Z}%
(n^{t}g_{n}(t),n)-\mathbf{Z}(n^{t_{k}}g_{n}(t_{k}),n)\right\Vert \right)
\neq 0\right) \\
&&\qquad \leq \mathbf{P}_{n}(\exists t,s\in \lbrack \gamma _{i},\gamma
_{i+1}-\varepsilon ]:\mathbf{Z}(n^{t}g_{n}(t),n)\neq \mathbf{Z}%
(n^{s}g_{n}(s),n)\,)
\end{eqnarray*}%
and to use the same arguments as before.

Theorem \ref{T_SkorohConst} is proved.

\subsection{Microscopic view\label{Sec52}}

We follow in this section the ideas of paper \cite{FZ} and to this aim
formulate a particular and slightly modified case of Theorem 6.5.4 in \cite%
{GS71} giving a convergence criterion in Skorokhod topology for a class of
Markov processes.

Let $\mathbf{K}_{n}(y),n=1,2,...$ be a sequence of Markov processes with
values in $\mathbb{Z}_{+}^{N}$ whose trajectories belong with probability 1
to the space $D_{[a,b]}(\mathbb{Z}_{+}^{N})$ of cadlag functions on $[a,b]$.

\begin{theorem}
\label{T_skoroh}If the finite-dimensional distributions of $\left\{ \mathbf{K%
}_{n}(y),a\leq y\leq b\right\} $ converge, as $n\rightarrow \infty ,$ to the
respective finite-dimensional distributions of a process $\left\{ \mathbf{K}%
(y),a\leq y\leq b\right\} $ and there exists a partition $\mathbb{Z}_{+}^{N}=%
\mathcal{B}\cup \mathcal{C},\mathcal{B}\cap \mathcal{C}=\varnothing $ such
that
\begin{equation*}
\lim_{h\downarrow 0}\overline{\lim_{n\rightarrow \infty }}\sup_{0\leq
s-y\leq h}\sup_{\mathbf{z}\in \mathcal{B}}\mathbf{P}(\mathbf{K}_{n}(s)\neq
\mathbf{K}_{n}(y)|\mathbf{K}_{n}(y)=\mathbf{z})=0,
\end{equation*}%
and
\begin{equation*}
\lim_{n\rightarrow \infty }\mathbf{P}(\exists y\in \lbrack a,b]:\mathbf{K}%
_{n}(y)\in \mathcal{C})=0
\end{equation*}%
then, as $n\rightarrow \infty $
\begin{equation*}
\mathcal{L}\left\{ \mathbf{K}_{n}(y),a\leq y\leq b\right\} \Longrightarrow
\mathcal{L}\left\{ \mathbf{K}(y),a\leq y\leq b\right\} .
\end{equation*}
\end{theorem}

In view of Lemma \ref{L_convol} the law\textbf{\ }$\mathbf{P}_{n}(\left\{
\mathbf{Z}(m,n),0\leq m\leq n\right\} \in (\cdot )|\mathbf{Z}(n)\neq \mathbf{%
0})$ specifies, for each fixed $n$ an inhomogeneous Markov branching
process. We denote its transition probabilities by $\mathbf{P}_{n}(m_{1},%
\mathbf{z};m_{2},(\cdot ))$.

Proving the tightness of $\mathbf{U}_{i}\left( \cdot \right) ,$ $%
i=1,2,...,N, $ we need to construct an appropriate partition of $\mathbb{Z}%
_{+}^{N}$ and to use Theorem \ref{T_skoroh} for each $\left[ 0,b\right]
\subset \lbrack 0,\infty )$.

Observe that if $\mathbf{w}=( w_{1},...,w_{N}) \leq \mathbf{z}=(
z_{1},...,z_{N}) $ (where the inequality is understood componentwise) then%
\begin{equation*}
\mathbf{P}_{n}( m_{0},\mathbf{w};m_{1},\left\{ \mathbf{w}\right\} ) \geq
\mathbf{P}_{n}( m_{0},\mathbf{z};m_{1},\left\{ \mathbf{z}\right\} ) .
\end{equation*}

Let $\ \mathcal{C}(k)=\left\{ \mathbf{z}\in \mathbb{Z}_{+}^{N}:\left\Vert
\mathbf{z}\right\Vert \leq k\right\} ,$%
\begin{equation*}
\mathcal{C}_{i}(k)=\left\{ \mathbf{z}\in \mathbb{Z}%
_{+}^{N}:z_{1}+...+z_{i-1}>0;\left\Vert \mathbf{z}\right\Vert \leq k\right\}
,~\mathcal{J}_{i}(k)=\mathcal{C}(k)\backslash \mathcal{C}_{i}(k).
\end{equation*}

Fix $i\in \{1,...,N-1\}$ and denote $m_{j}=( Y_{j}+l_{n}) n^{\gamma
_{i}},j=1,2$.

\begin{lemma}
\label{L_ceretain event.}Under Hypothesis A for any fixed $k$ and $%
0<b<\infty $%
\begin{equation*}
\lim_{h\downarrow 0}\overline{\lim_{n\rightarrow \infty}}\sup_{\substack{ %
0\leq Y_{1}-Y_{0}\leq h,  \\ Y_{1},Y_{0}\in \lbrack 0,b]}}\sup_{\mathbf{z}%
\in \mathcal{J}_{i}(k)}\mathbf{P}_{n}(\mathbf{Z}(m_{1};n)\neq \mathbf{z}|%
\mathbf{Z}(m_{0};n)=\mathbf{z})=0.
\end{equation*}
\end{lemma}

\textbf{Proof.} By the branching property, the decomposability assumption,
and the positivity of the offspring number of each particle in the reduced
process we have for all $m_{1}\geq m_{0}$ and $\mathbf{z}\in \mathcal{J}%
_{i}(k)$
\begin{eqnarray*}
\mathbf{P}_{n}(m_{0},\mathbf{z};m_{1},\left\{ \mathbf{z}\right\} )
&=&\prod_{j=i}^{N}(\mathbf{P}_{n}(m_{0},\mathbf{e}_{j};m_{1},\left\{ \mathbf{%
e}_{j}\right\} ))^{z_{j}} \\
&\geq &\prod_{j=i}^{N}(\mathbf{P}_{n}(m_{0},\mathbf{e}_{j};m_{1},\left\{
\mathbf{e}_{j}\right\} ))^{k}.
\end{eqnarray*}%
Using Lemma \ref{L_derivative} we get for $m_{0}=(Y_{0}+l_{n})n^{\gamma
_{i}} $ and $m_{1}=(Y_{1}+l_{n})n^{\gamma _{i}}$:
\begin{equation}
\inf_{\substack{ 0\leq Y_{1}-Y_{0}\leq h,  \\ Y_{0},Y_{1}\in \left[ 0,b%
\right] }}\mathbf{P}_{n}(m_{0},\mathbf{z};m_{1},\left\{ \mathbf{z}\right\}
)\geq (1-\chi h)^{Nk}.  \label{kapBel}
\end{equation}%
This implies the claim of the lemma. \

\begin{lemma}
\label{L_ratio}If $m_{j}=(Y_{j}+l_{n})n^{\gamma _{i}},~j=0,1,2,$ and $0\leq
Y_{0}<Y_{1}<Y_{2}$ with $Y_{1}-Y_{0}\leq h$, then for all $n\geq n_{0}$
\begin{eqnarray*}
&&\mathbf{P}_{n}(\mathbf{Z}(m_{1},n)=\mathbf{z}|\mathbf{Z}(m_{0},n)=\mathbf{z%
};\left\Vert \mathbf{Z}(m_{2},n)\right\Vert \leq k) \\
&&\qquad \geq \mathbf{P}_{n}(m_{0},\mathbf{z};m_{1},\left\{ \mathbf{z}%
\right\} )\frac{\mathbf{P}_{n}(m_{1},\mathbf{z};m_{2},\mathcal{C}(k))}{%
\mathbf{P}_{n}(m_{1},\mathbf{z};m_{2},\mathcal{C}(k))\mathbf{+}\chi Nkh}.
\end{eqnarray*}
\end{lemma}

\textbf{Proof.} We have%
\begin{eqnarray*}
&&\mathbf{P}_{n}(\mathbf{Z}(m_{1},n)=\mathbf{z}|\mathbf{Z}(m_{0},n)=\mathbf{z%
},\left\Vert \mathbf{Z}(m_{2},n)\right\Vert \leq k) \\
&&\qquad =\mathbf{P}_{n}(m_{0},\mathbf{z};m_{1},\left\{ \mathbf{z}\right\} )%
\frac{\mathbf{P}_{n}(m_{1},\mathbf{z};m_{2},\mathcal{C}(k))}{\mathbf{P}%
_{n}(m_{0},\mathbf{z};m_{2},\mathcal{C}(k))}.
\end{eqnarray*}%
In view of (\ref{kapBel})%
\begin{eqnarray*}
\mathbf{P}_{n}(m_{0},\mathbf{z};m_{2},\mathcal{C}(k)) &=&\sum_{\mathbf{z}%
^{\prime }}\mathbf{P}_{n}(m_{0},\mathbf{z};m_{1},\left\{ \mathbf{z}^{\prime
}\right\} )\mathbf{P}_{n}(m_{1},\mathbf{z}^{\prime };m_{2},\mathcal{C}(k)) \\
&\leq &1-\mathbf{P}_{n}(m_{0},\mathbf{z};m_{1},\left\{ \mathbf{z}\right\} )+%
\mathbf{P}_{n}(m_{1},\mathbf{z};m_{2},\mathcal{C}(k)) \\
&\leq &1-(1-\chi h)^{Nk}+\mathbf{P}_{n}(m_{1},\mathbf{z};m_{2},\mathcal{C}%
(k)) \\
&\leq &\chi Nkh+\mathbf{P}_{n}(m_{1},\mathbf{z};m_{2},\mathcal{C}(k)).
\end{eqnarray*}%
Hence the needed statement follows.

\begin{lemma}
\label{L_Kskrokh}Under Hypothesis A for any fixed $k,$ $0<b<\infty ,$ and $%
m_{0}=( Y_{0}+l_{n}) n^{\gamma _{i}},\ m_{1}=( Y_{1}+l_{n}) n^{\gamma
_{i}},\ m_{2}=2bn^{\gamma _{i}}$ we have%
\begin{equation*}
\lim_{h\downarrow 0}\overline{\lim_{n\rightarrow \infty }}\sup_{\substack{ %
0\leq Y_{1}-Y_{0}\leq h  \\ Y_{1},Y_{0}\in \lbrack 0,b]}}\sup_{\mathbf{z}\in
\mathcal{J}_{i}(k)}\mathbf{P}_{n}( \mathbf{Z}( m_{1};n) \neq \mathbf{z}|%
\mathbf{Z}( m_{0};n) =\mathbf{z},\left\Vert \mathbf{Z}( m_{2},n) \right\Vert
\leq k) =0.
\end{equation*}
\end{lemma}

\textbf{Proof.} By (\ref{kapBel}) and Lemma \ref{L_ratio} for $%
m_{0}=(Y_{0}+l_{n})n^{\gamma _{i}}$ and $m_{1}=(Y_{1}+l_{n})n^{\gamma _{i}}$
\begin{eqnarray*}
&&\mathbf{P}_{n}(\mathbf{Z}(m_{1},n)=\mathbf{z}|\mathbf{Z}(m_{0},n)=\mathbf{z%
};\left\Vert \mathbf{Z}(m_{2},n)\right\Vert \leq k) \\
&&\qquad \geq (1-\chi h)^{Nk}\frac{\mathbf{P}_{n}(m_{1},\mathbf{z};m_{2},%
\mathcal{C}(k))}{\mathbf{P}_{n}(m_{1},\mathbf{z};m_{2},\mathcal{C}(k))%
\mathbf{+}\chi Nkh}.
\end{eqnarray*}%
Using the decomposability hypothesis and Lemma \ref{L_derivative} we obtain
\begin{eqnarray*}
&&\mathbf{P}_{n}(m_{1},\mathbf{z};m_{2},\mathcal{C}(k))\geq \mathbf{P}%
_{n}(m_{1},\mathbf{z};m_{2},\left\{ \mathbf{z}\right\} ) \\
&&\qquad =\prod_{j=i}^{N}(\mathbf{P}_{n}(m_{1},\mathbf{e}_{j};m_{2},\left\{
\mathbf{e}_{j}\right\} ))^{z_{j}}\geq \prod_{j=i}^{N}\mathbf{P}%
_{n}^{k}(l_{n}n^{\gamma _{i}},\mathbf{e}_{j};2bn^{\gamma _{i}},\left\{
\mathbf{e}_{j}\right\} ).
\end{eqnarray*}%
It follows from Theorem~\ref{T_findim} that
\begin{equation*}
\lim_{n\rightarrow \infty }\prod_{j=i}^{N}\mathbf{P}_{n}^{k}(l_{n}n^{\gamma
_{i}},\mathbf{e}_{j};2bn^{\gamma _{i}},\left\{ \mathbf{e}_{j}\right\} )=%
\mathbf{P}^{k}(\mathbf{U}_{i}(2b)=\mathbf{e}_{i}|\mathbf{U}_{i}(0)=\mathbf{e}%
_{i})=B>0.
\end{equation*}%
Hence we get%
\begin{eqnarray*}
&&\varliminf_{n\rightarrow \infty }\inf_{\substack{ 0\leq Y_{1}-Y_{0}\leq h
\\ Y_{1},Y_{0}\in \lbrack 0,b]}}\inf_{\mathbf{z}\in \mathcal{J}_{i}(k)}%
\mathbf{P}_{n}(\mathbf{Z}(m_{1};n)=\mathbf{z}|\mathbf{Z}(m_{0};n)=\mathbf{z}%
,\left\Vert \mathbf{Z}(m_{2},n)\right\Vert \leq k) \\
&&\qquad \qquad \geq (1-\chi h)^{Nk}\frac{B}{B\mathbf{+}\chi Nkh}.
\end{eqnarray*}%
Letting $h\downarrow 0$ completes the proof of the lemma. \

\begin{corollary}
\label{C_skk}Under the conditions of Lemma \ref{L_Kskrokh}
\begin{eqnarray*}
&&\mathcal{L}\left\{ \mathbf{Z}((y+l_{n})n^{1/2^{N-i}},n),0\leq y\leq b\,%
\Big|\,\left\Vert \mathbf{Z}(m_{2},n)\right\Vert \leq k,\mathbf{Z}(n)\neq
\mathbf{0}\right\} \\
&&\qquad\qquad\qquad\qquad\qquad\quad\,\Longrightarrow \mathcal{L}%
_{R_{i}}\left\{ \mathbf{U}_{i}(y),0\leq y\leq b|\,\left\Vert \mathbf{U}%
_{i}(2b)\right\Vert \leq k\right\} .
\end{eqnarray*}
\end{corollary}

\textbf{Proof.} Convergence of finite-dimensional distributions follows from
the respective results for the convergence of the processes established in
point 1) of Theorem \ref{T_Skhod1}. Tightness follows from Lemma \ref%
{L_Kskrokh} and Theorem \ref{T_skoroh} by taking $\mathcal{B}=\mathcal{J}%
_{i}(k)$ and $\mathcal{C=C}_{i}(k)$ and observing that
\begin{eqnarray*}
&&\lim_{n\rightarrow \infty }\mathbf{P}_{n}(\mathbf{Z}(l_{n}n^{\gamma
_{i}},n)\in \mathcal{C}_{i}(k)|\left\Vert \mathbf{Z}(m_{2},n)\right\Vert
\leq k) \\
&&\qquad \leq \lim_{n\rightarrow \infty }\mathbf{P}_{n}\left( \text{\b{Z}}%
_{i-1}(l_{n}n^{\gamma _{i}})>0|\left\Vert \mathbf{Z}(m_{2},n)\right\Vert
\leq k\right) =0.
\end{eqnarray*}

\textbf{Proof} \textbf{of Theorem \ref{T_Skhod1}.} Let for $c>b$
\begin{eqnarray*}
\mathbf{P}_{n,i}(b;(\cdot )) &=&\mathbf{P}_{n}(\left\{ \mathbf{Z}%
((y+l_{n})n^{\gamma _{i}},n),0\leq y\leq b\right\} \in (\cdot )), \\
\mathbf{P}_{n,i}^{(k)}(b,c;(\cdot )) &=&\mathbf{P}_{n}(\left\{ \mathbf{Z}%
((y+l_{n})n^{\gamma _{i}},n),0\leq y\leq b\right\} \in (\cdot )|\left\Vert
\mathbf{Z}(cn^{\gamma _{i}},n)\right\Vert \leq k), \\
\mathbf{\bar{P}}_{n,i}^{(k)}(b,c;(\cdot )) &=&\mathbf{P}_{n}(\left\{ \mathbf{%
Z}((y+l_{n})n^{\gamma _{i}},n),0\leq y\leq b\right\} \in (\cdot )|\left\Vert
\mathbf{Z}(cn^{\gamma _{i}},n)\right\Vert >k)
\end{eqnarray*}%
and%
\begin{eqnarray*}
\mathcal{P}_{i}(b;(\cdot )) &=&\mathbf{P}_{R_{i}}(\left\{ \mathbf{U}%
_{i}(y),0\leq y\leq b\right\} \in (\cdot )), \\
\mathcal{P}_{i}^{(k)}(b,c;(\cdot )) &=&\mathbf{P}_{R_{i}}(\left\{ \mathbf{U}%
_{i}(y),0\leq y\leq b\right\} \in (\cdot )|\left\Vert \mathbf{U}%
_{i}(c)\right\Vert \leq k).
\end{eqnarray*}%
Then for $0<b<\infty $ and a continuous real function $\psi $ on $D_{[0,b]}(%
\mathbb{Z}_{+}^{N})$ such that $\left\vert \psi \right\vert \leq q$ for a
positive $q$ we have
\begin{eqnarray*}
\int \psi (x)\mathbf{P}_{n,i}(b;dx) &=&\mathbf{P}_{n}(\left\Vert \mathbf{Z}%
(2bn^{\gamma _{i}},n)\right\Vert >k)\int \psi (x)\mathbf{\bar{P}}%
_{n,i}^{(k)}(b,2b;dx) \\
&&+\mathbf{P}_{n}(\left\Vert \mathbf{Z}(2bn^{\gamma _{i}},n)\right\Vert \leq
k)\int \psi (x)\mathbf{P}_{n,i}^{(k)}(b,2b;dx).
\end{eqnarray*}%
For the first summand we get%
\begin{eqnarray*}
&&\lim \sup_{n\rightarrow \infty }\mathbf{P}_{n}(\left\Vert \mathbf{Z}%
(2bn^{\gamma _{i}},n)\right\Vert >k)\int \psi (x)\mathbf{\bar{P}}%
_{n,i}^{(k)}(b,2b;dx) \\
&&\quad \leq q\lim \sup_{n\rightarrow \infty }\mathbf{P}_{n}(\left\Vert
\mathbf{Z}(2bn^{\gamma _{i}},n)\right\Vert >k)=q\mathbf{P}%
_{R_{i}}(\left\Vert \mathbf{U}_{i}(2b)\right\Vert >k)=o(1)
\end{eqnarray*}%
as $k\rightarrow \infty $ by the properties of $\mathbf{U}_{i}(\cdot )$. \

On the other hand, letting first $n\rightarrow \infty $ and than $%
k\rightarrow \infty $ we obtain%
\begin{eqnarray*}
&&\lim_{k\rightarrow \infty }\lim_{n\rightarrow \infty }\mathbf{P}%
_{n}(\left\Vert \mathbf{Z}(2bn^{\gamma _{i}},n)\right\Vert \leq k)\int \psi
(x)\mathbf{P}_{n,i}^{(k)}(b,2b;dx) \\
&&\quad =\lim_{k\rightarrow \infty }\mathbf{P}_{R_{i}}(0<\left\Vert \mathbf{U%
}_{i}(2b)\right\Vert \leq k)\int \psi (x)\mathcal{P}_{i}^{(k)}(b,2b;dx) \\
&&\quad =\lim_{k\rightarrow \infty }\int_{\left\{ 0<\left\Vert \mathbf{U}%
_{i}(2b)\right\Vert \leq k\right\} }\psi (x)\mathcal{P}_{i}(b,2b;dx)=\int
\psi (x)\mathcal{P}_{i}(b;dx).
\end{eqnarray*}%
Thus,%
\begin{equation*}
\lim_{n\rightarrow \infty }\int \psi (x)\mathbf{P}_{n}(\mathbf{Z}((\cdot
+l_{n})n^{\gamma _{i}},n)\in dx)=\int \psi (x)\mathcal{P}_{i}(b;dx)
\end{equation*}%
for any bounded continuous function on $D_{[0,b]}(\mathbb{Z}_{+}^{N})$
proving point 1) of Theorem~\ref{T_Skhod1}.

The proof of point 2) of Theorem \ref{T_Skhod1} needs only a few changes in
comparison with the proof of the respective theorem in \cite{FZ} and we omit
it. \

\section{Proofs of Theorems \protect\ref{T_mrcaMany} and \protect\ref{T_type}
\label{Sec7}}

\textbf{Proof} \textbf{of Theorem \ref{T_mrcaMany}.} Our arguments are based
on the following simple observation
\begin{equation*}
\left\{ \bar{Z}_{1}(m,n)=1\right\} \Leftrightarrow \left\{ \beta _{n}\geq
m\right\} .
\end{equation*}

\textbf{Proof of 1).} According to (\ref{LimSing}) for $m\ll n^{\gamma _{1}}$
\begin{eqnarray*}
&&\lim_{n\rightarrow \infty }\mathbf{P}_{n}(\bar{Z}_{1}(m,n)=1)=\lim_{n%
\rightarrow \infty }\mathbf{P}_{n}(Z_{1}(m,n)=1) \\
&&\quad +\lim_{n\rightarrow \infty }\mathbf{P}_{n}(\bar{Z}_{2}(m,n)=1)=1+0=1.
\end{eqnarray*}

\textbf{Proof of 2).} Observe that by point 2) of Theorem \ref{T_familyMany}%
\begin{eqnarray*}
&&\lim_{n\rightarrow \infty }\mathbf{P}_{n}(\beta _{n}\geq yn^{\gamma
_{i}})=\lim_{n\rightarrow \infty }\mathbf{P}_{n}(\bar{Z}_{1}(yn^{\gamma
_{i}},n)=1) \\
&&\quad =\lim_{n\rightarrow \infty }\mathbf{P}_{n}(Z_{i}(yn^{\gamma
_{i}},n)+Z_{i+1}(yn^{\gamma _{i}},n)=1) \\
&&\quad =\lim_{n\rightarrow \infty }\mathbf{P}_{n}(Z_{i}(yn^{\gamma
_{i}},n)=1)+\lim_{n\rightarrow \infty }\mathbf{P}_{n}(Z_{i+1}(yn^{\gamma
_{i}},n)=1).
\end{eqnarray*}%
Direct calculations show that%
\begin{equation*}
-\frac{\partial \varphi _{i}(y;s_{i},s_{i+1})}{\partial s_{i}}\left\vert
_{s_{i}=s_{i+1}=0}\right. =\frac{1-\tanh (yb_{i}c_{iN})}{1+\tanh
(yb_{i}c_{iN})}=e^{-2yb_{i}c_{iN}}
\end{equation*}%
and%
\begin{equation*}
-\frac{\partial \varphi _{i}(y;s_{i},s_{i+1})}{\partial s_{i+1}}\left\vert
_{s_{i}=s_{i+1}=0}\right. =\frac{\tanh (yb_{i}c_{iN})}{1+\tanh (yb_{i}c_{iN})%
}=\frac{1}{2}-\frac{1}{2}\,e^{-2yb_{i}c_{iN}}.
\end{equation*}%
Thus,
\begin{eqnarray*}
&&\lim_{n\rightarrow \infty }\mathbf{P}_{n}(Z_{i}(yn^{\gamma
_{i}},n)=1;\beta _{n}\geq yn^{\gamma _{i}}) \\
&&\qquad \qquad =-\frac{\partial (\varphi _{i}(y;s_{i},s_{i+1}))^{1/2^{i-1}}%
}{\partial s_{i}}\left\vert _{s_{i}=s_{i+1}=0}\right. =\frac{1}{2^{i-1}}%
\,e^{-2yb_{i}c_{iN}}
\end{eqnarray*}%
and%
\begin{eqnarray*}
&&\lim_{n\rightarrow \infty }\mathbf{P}_{n}(Z_{i+1}(yn^{\gamma
_{i}},n)=1;\beta _{n}\geq yn^{\gamma _{i}}) \\
&&\qquad \quad =-\frac{\partial (\varphi _{i}(y;s_{i},s_{i+1}))^{1/2^{i-1}}}{%
\partial s_{i+1}}\left\vert _{s_{i}=s_{i+1}=0}\right. =\frac{1}{2^{i}}%
(1-e^{-2yb_{i}c_{iN}}).
\end{eqnarray*}%
Combining the previous estimates yields
\begin{equation*}
\lim_{n\rightarrow \infty }\mathbf{P}_{n}(\beta _{n}\leq yn^{\gamma _{i}})=1-%
\frac{1}{2^{i}}-\frac{1}{2^{i}}e^{-2yb_{i}c_{iN}}.
\end{equation*}

\textbf{Proof of 3).} This is evident.

\textbf{Proof of 4).} The needed statement follows from the equality%
\begin{equation*}
-\frac{\partial }{\partial s_{N}}\left( \frac{1}{x+(1-x)/(1-s_{N})}\right)
^{2^{-(N-1)}}\left\vert _{s_{N}=0}\right. =\frac{1}{2^{N-1}}(1-x).
\end{equation*}

\textbf{Proof} \textbf{of Theorem \ref{T_type}.} Consider the case $N\geq 4$
and $i\in \left\{ 2,3,..,N-2\right\} $ only. For $N=2,3$ or $N\geq 4$ and $%
i\in \left\{ 1,N-1\right\} $ some of the random variables (events) below do
not exist (are empty) and the needed arguments become shorter.

Since the total number of particles of all types in the reduced process does
not decrease with time, $\mathbf{P}_{n}(\beta _{n}<m)=\mathbf{P}_{n}\left(
\text{\b{Z}}_{1}(m,n)\geq 2\right) $. We now take%
\begin{equation*}
m_{i}=n^{\gamma _{i}(1+\gamma _{i})},\ i=1,2,...,N-1,
\end{equation*}%
and denote $\mathcal{H}_{i}=\left\{ m:m_{i-1}\leq m\leq m_{i}\right\} $.

Since $\bar{Z}_{i}(k,n)$ is monotone increasing in $k$ for each fixed $n,$
Theorem \ref{T_familyMany} and~(\ref{SurvivSingle}) imply, as $n\rightarrow
\infty $
\begin{eqnarray*}
&&\mathbf{P}_{n}(\zeta _{n}=i;\beta _{n}\notin \mathcal{H}_{i})\leq \mathbf{P%
}_{n}(\exists k<m_{i-1}:Z_{i}(k,n)>0) \\
&&\qquad +\mathbf{P}_{n}(\exists k>m_{i}:Z_{i}(k,n)>0) \\
&&\quad \leq \mathbf{P}_{n}(\bar{Z}_{i}(m_{i-1},n)>0)+\mathbf{P}%
_{n}(Z_{i}(m_{i})>0)=o(1).
\end{eqnarray*}%
By the same statements we conclude, as $n\rightarrow \infty $%
\begin{eqnarray*}
&&\mathbf{P}_{n}(\zeta _{n}\notin \{i,i+1\};\beta _{n}\in \mathcal{H}%
_{i})\leq \mathbf{P}_{n}\left( \exists k\in \mathcal{H}_{i}:\text{\b{Z}}%
_{i-1}(k,n)+\bar{Z}_{i+2}(k,n)>0\right) \\
&&\quad \leq \mathbf{P}_{n}\left( \exists k\in \mathcal{H}_{i}:\text{\b{Z}}%
_{i-1}(k)+\bar{Z}_{i+2}(k,n)>0\right) \\
&&\quad \leq \mathbf{P}_{n}(\exists k\in \mathcal{H}_{i}:\text{\b{Z}}%
_{i-1}(k)>0)+\mathbf{P}_{n}\left( \exists k\in \mathcal{H}_{i}:\bar{Z}%
_{i+2}(k,n)>0\right) \\
&&\quad \leq \mathbf{P}_{n}(\text{\b{Z}}_{i-1}(m_{i-1})>0)+\mathbf{P}_{n}(%
\bar{Z}_{i+2}(m_{i},n)>0)=o(1).
\end{eqnarray*}%
Hence, as $n\rightarrow \infty $%
\begin{eqnarray}
\mathbf{P}_{n}(\zeta _{n}=i) &=&\mathbf{P}_{n}(\zeta _{n}=i;\beta _{n}\in
\mathcal{H}_{i})+o(1)  \notag \\
&=&\mathbf{P}_{n}(\beta _{n}\in \mathcal{H}_{i})-\mathbf{P}_{n}(\zeta
_{n}=i+1;\beta _{n}\in \mathcal{H}_{i})+o(1).  \label{Type}
\end{eqnarray}%
Introduce the event
\begin{equation*}
\mathcal{G}_{i}(j,n)=\left\{ \text{\b{Z}}_{i}(j;n)+\bar{Z}%
_{i+2}(j+1,n)=0;Z_{i+1}(j,n)=1\right\} .
\end{equation*}%
Clearly,%
\begin{eqnarray*}
&&\mathbf{P}_{n}(\zeta _{n}=i+1;\beta _{n}\in \mathcal{H}_{i})=%
\sum_{j=m_{i-1}}^{m_{i}}\mathbf{P}_{n}(\zeta _{n}=i+1;\beta _{n}=j) \\
&&\,=\sum_{j=m_{i-1}}^{m_{i}}\mathbf{P}_{n}(\mathcal{G}_{i}(j,n),\bar{Z}%
_{i+1}(j+1,n)\geq 2) \\
&&\,=o(1)+\sum_{j=m_{i-1}}^{m_{i}}\mathbf{P}_{n}(\mathcal{G}_{i}(j,n))%
\mathbf{P}_{n}(Z_{i+1}(j+1,n)\geq 2|\mathbf{Z}(j,n)=\mathbf{e}_{i+1}).
\end{eqnarray*}%
It is not difficult to check (recall (\ref{DefNONimmigr}), (\ref{Derivat})
and (\ref{ExplCoeff})) that
\begin{eqnarray*}
\mathbf{P}_{n}(Z_{i+1}(j+1,n)=1|\mathbf{Z}(j,n)=\mathbf{e}_{i+1}) &=&\frac{%
Q_{n-j-1}^{(i+1,N)}}{Q_{n-j}^{(i+1,N)}}\frac{dh_{i+1}(s,\mathbf{1}^{(N-i-1)})%
}{ds}\left\vert _{s=H_{n-j-1}^{(i+1,N)}(\mathbf{0})}\right. \\
&\geq &\frac{dh_{i+1}(s,\mathbf{1}^{(N-i-1)})}{ds}\left\vert
_{s=H_{n-j-1}^{(i+1,N)}(\mathbf{0})}\right. \\
&\geq &1-2b_{i+1}Q_{n-j-1}^{(i+1,N)} \\
&\geq &1-2b_{i+1}Q_{n-m_{i}}^{(i+1,N)}.
\end{eqnarray*}%
Hence, using the estimate%
\begin{eqnarray*}
\mathbf{P}_{n}(Z_{i+1}(j+1,n)\geq 2|\mathbf{Z}(j,n)=\mathbf{e}_{i+1}) &=&1-%
\mathbf{P}_{n}(Z_{i+1}(j+1,n)=1|\mathbf{Z}(j,n)=\mathbf{e}_{i+1}) \\
&\leq &2b_{i}Q_{n-m_{i}}^{(i+1,N)}
\end{eqnarray*}%
we conclude%
\begin{eqnarray*}
\mathbf{P}_{n}(\zeta _{n}=i+1;\beta _{n}\in \mathcal{H}_{i})
&=&o(1)+O(m_{i}Q_{n-m_{i}}^{(i+1,N)}) \\
&=&o(1)+O(n^{\gamma _{i}(1+\gamma _{i})}n^{-\gamma _{i+1}})=o(1).
\end{eqnarray*}%
This, on account of (\ref{recent_i}) and (\ref{Type}) gives
\begin{equation*}
\lim_{n\rightarrow \infty }\mathbf{P}_{n}(\zeta _{n}=i)=\lim_{n\rightarrow
\infty }\mathbf{P}_{n}(\beta _{n}\in \mathcal{H}_{i})=\lim_{n\rightarrow
\infty }\mathbf{P}_{n}(n^{\gamma _{i}}\ll \beta _{n}\ll n^{\gamma _{i+1}})=%
\frac{1}{2^{i}}
\end{equation*}%
as desired.

Finally,%
\begin{equation*}
\lim_{n\rightarrow \infty }\mathbf{P}_{n}( \zeta _{n}=N) =1-\sum_{i=1}^{N-1}%
\frac{1}{2^{i}}=\frac{1}{2^{N-1}}.
\end{equation*}

Theorem \ref{T_type} is proved.

\textbf{Acknowledgement}. This work was partially supported by the Russian
Foundation for Basic Research, project N14-01-00318. The author would also
like to thank prof. A.M.Zubkov for valuable remarks.

\end{document}